\documentclass[11pt,reqno]{amsart}
\usepackage{textcomp}
\usepackage{amssymb, a4wide}
\usepackage{amsmath}
\usepackage{graphicx}
\usepackage{epsfig}
\usepackage{subcaption}
\DeclareMathOperator{\arcsinh}{arcsinh}
\DeclareMathOperator{\arctanh}{arctanh}
\DeclareMathOperator{\arccosh}{arccosh}

\newtheorem{theorem}{Theorem}
\newtheorem{example}{Example}
\newtheorem{remark}{Remark}
\newtheorem{lemma}{Lemma}

\newtheorem{definition}{Definition}
\usepackage{caption}
\usepackage{subcaption}

\title{Isometric Timelike Surfaces in 4--Dimensional Minkowski Space}

\author[B. Bekta\c{s} Dem\.{i}rc\.i]{Burcu Bekta\c s Dem\.{i}rc\.i}
\address{Fatih Sultan Mehmet Vak{\i}f University, Hal\.{I}\c{c} Campus, Faculty of Engineering,
Department of Software Engineering, 34445, Beyo\u{g}lu, \.{I}stanbul, T{\"u}rk\.{I}ye}
\email{bbektas@fsm.edu.tr}

\author[M. Babaarslan]{Murat Babaarslan}
\address{Yozgat Bozok University, Department of Mathematics, 66100, Yozgat, T{\"u}rk\.{I}ye}
\email{murat.babaarslan@bozok.edu.tr}

\author[Y. K{\"u}\c{c}{\"u}kar{\i}kan]{Yas\.in K{\"u}\c{c}{\"u}kar{\i}kan}
\address{Yozgat Bozok University, School of Graduates Studies, Department of Mathematics, 66100, Yozgat, T{\"u}rk\.{I}ye}
\email{kucukarikanyasin@gmail.com}

\begin{document}

\begin{abstract}
In this paper, first we study on Bour's theorem for four kinds of timelike helicoidal surfaces in 4-dimensional Minkowski space. Secondly, we analyse the geometric properties of these isometric surfaces having same Gauss map. Also, we present the parametrizations of such isometric pair of surfaces. Finally, we introduce some examples and draw the corresponding graphs by using Wolfram Mathematica 10.4. 
\end{abstract}

\maketitle
\section{Introduction}

One of the most important knowledge in the surface theory is that the right helicoid and catenoid is only minimal ruled surface and minimal rotational surface, respectively. Also, it is known that they have same Gauss map \cite{Ikawa1}. In the surface theory, the following Bour's theorem is quite popular:

\textbf{Bour's theorem.}\cite{Bour}
A generalized helicoid is isometric to a rotational surface so that helices on the helicoid correspond to parallel circles on the rotational surface. 

In 2000, Ikawa \cite{Ikawa1} gave the parametrizations of the pairs of surface of Bour's theorem which have same Gauss map in 3-dimensional Euclidean space $\mathbb{E}^3$.  Helicoidal surfaces with constant mean curvature in $\mathbb{E}^3$ were investigated by do Carmo and Dajczer \cite{Carmo}. Also, spacelike helicoidal surfaces with constant mean curvature in 3-dimensional Minkowski space $\mathbb{E}^3_1$ were studied by Sasahara \cite{Sasahara}. In 2002, Ikawa \cite{Ikawa2} studied Bour's theorem for spacelike and timelike generalized helicoid  with non--null and null axis in $\mathbb{E}^3_1$. 
Bour's theorem for generalized helicoid with null axis in $\mathbb{E}^3_1$ was introduced by G{\"u}ler and Vanl{\i} \cite{Guler1} in 2006. In 2010, G{\"u}ler et al. \cite{Guler3} investigated Bour's theorem for the Gauss map of generalized helicoid in $\mathbb{E}^3$. As a generalization, in 2015, Bour's theorem for helicoidal surfaces in $\mathbb{E}^3$ were studied by G{\"u}ler and Yayl{\i} \cite{Guler3}.

In 2017, Hieu and Thang \cite{Hieu} studied on Bour's theorem for helicoidal surfaces 
in 4-dimensional Euclidean space $\mathbb{E}^4$ and 
they proved that if the Gauss maps of isometric surfaces are same, then they are hyperplanar and minimal. 
Also, they gave the parametrizations of such minimal surfaces.

Einstein's theory of special relativity is strongly related to Minkowski space--time (or called as 4-dimensional Minkowski space) $\mathbb{E}^4_1$ (see \cite{Ratcliffe} for details). Because of this important relation, in 2021, 
Babaarslan and Sönmez \cite{Babaarslan1} gave the parametrizations of three types of helicoidal surfaces in $\mathbb{E}^4_1$ 
by using three types rotation with $2$--dimensional axis, called elliptic, hyperbolic and parabolic rotations which leave the spacelike, timelike and lightlike planes invariant. Bour's theorem for these spacelike helicodial surfaces in $\mathbb{E}^4_1$ were introduced by Babaarslan et al. \cite{Babaarslan2}.

In this paper, we continue to study on Bour's theorem for four kinds of timelike helicoidal surfaces in $\mathbb{E}^4_1$. We analyse the geometric properties of these timelike isometric surfaces having same Gauss map as hyperplanar and minimal. Also, we give the parametrizations of such timelike isometric pair of surfaces. Finally, we give some examples by using Wolfram Mathematica 10.4.

\section{Preliminaries}
In this subsection, we recall some basic definitions and formulas in 4-dimensional Minkowski space $\mathbb{E}^4_1$. For more information, we refer to \cite{ONeill}.

A metric tensor $g$ is symmetric, bilinear, non-degenarate and (0,2) tensor field in $\mathbb{E}^4_1$ which is defined by
\begin{equation}
\label{eq1}
g(x,y)=\langle x, y\rangle=x_{1}y_{1}+x_{2}y_{2}+x_{3}y_{3}-x_{4}y_{4}
\end{equation}
for the vectors  $x=(x_{1},x_{2},x_{3},x_{4}), y=(y_{1},y_{2},y_{3},y_{4}) \in \mathbb{E}^4_{1}$.

The causal character of a vector $x \in \mathbb{E}^4_1$ is spacelike if $\langle x,x\rangle>0$ or $x=0$, timelike if $\langle x,x\rangle<0$ and lightlike (null) if $\langle x,x\rangle=0$ and $x\neq 0$.

A curve in $\mathbb{E}^4_{1}$ is a smooth mapping $\alpha:I\subset\mathbb{R}\longrightarrow \mathbb{E}^4_1$, where $I$ is an open interval. The tangent vector of $\alpha$ at $t \in I $ is given by $\alpha'(t)$ and $\alpha$ is a regular curve if $\alpha'(t)\neq 0$ for all $t$. Also, $\alpha$ is spacelike if all of its tangent vectors $\alpha'(t)$ spacelike; similarly for lightlike and timelike.

\begin{definition}\cite{Kaya}
\label{circledef}
We suppose that the plane $P$ involving the circle is the plane of equation, $x_3=0$, $x_1=0$ or $x_2-x_3=0$, if $P$ is spacelike, timelike or lightlike, respectively. Thus, a circle $C \in \mathbb{E}^3_1$ can be defined as follows:
\begin{itemize}
\item If $P \equiv \{x_3=0\}$, then $C$ is an Euclidean circle $\alpha(s)=p+ r(\cos{s},\sin{s},0)$ with center $p\in P$ and radius $r>0$.
\item If $P \equiv \{x_1=0\}$, then $C$ is a spacelike hyperbola $\alpha(s)=p+ r(0,\sinh{s},\cosh{s})$ or $C$ is a timelike hyperbola $\alpha(s)=p+ r(0,\cosh{s},\sinh{s})$, where $p\in P$ and $r>0$ is the radius.
\item If $P \equiv \{x_2-x_3=0\}$, then $C$ is spacelike parabola $\alpha(s)=p+ (s,rs^2,rs^2)$, where $p\in P$ and $r>0$.
\end{itemize}
\end{definition}

A semi-Riemann surface $X$ is a 2-dimensional semi-Riemann manifold in $\mathbb{E}^4_1$. 
For a coordinate system $\{u,v\}$ in $X$, the tangent plane of $X$ at $p$ is given by $T_{p}X=span\{X_{u}, X_{v}\}$. The components of the metric tensor are denoted by
\begin{equation}
\label{eq2}
g_{11}=\langle X_{u}, X_{u}\rangle,
\text{ }g_{12}=g_{21}=\langle X_{u}, X_{v}\rangle, 
\text{ }g_{22}=\langle X_{v}, X_{v}\rangle. 
\end{equation}
Thus, the first fundamental form (or line element) is
\begin{equation}
g=g_{11} d{u}^{2}+ 2g_{12}d{u}d{v}+g_{22}dv^2.
\end{equation}

When ${W}=\mbox{det}(g)=g_{11}g_{22}-
g_{12}^2\not=0 $, the semi-Riemann surface $X$ is non--degenerate, namely, when ${W}>0$, the semi-Riemann surface $X$ is spacelike and when ${W}<0$, $X$ is a timelike surface. 

Let $\{e_{1},{e}_{2},{N}_{1},{N}_{2}\}$ be 
a local orthonormal frame on the semi-Riemann surface $X$ in $\mathbb{E}^4_1$
such that 
$e_{1},{e}_{2}$ are tangent to $X$ and $N_{1},N_{2}$ are normal to $X$. 
The coefficients of the second fundamental form tensor according to $N_{i}$, $(i=1,2)$ are denoted by
 \begin{equation}
 \label{eq3}
 b_{11}^{i}=\langle X_{uu},  N_{i}\rangle,
 \;\;b_{12}^{i}=b_{21}^{i}=\langle X_{uv}, N_{i}\rangle, 
 \;\;b_{22}^{i}=\langle X_{vv}, N_{i}\rangle.  
 \end{equation}
The mean curvature vector $H$ of $X$ in $\mathbb{E}^4_1$ is given by
\begin{equation}
\label{eq4}
H={\epsilon}_{1}H_{1}N_{1}+\mathbf{\epsilon}_{2}H_{2}N_{2},
\end{equation}
where the components $H_i$ of $H$ is
$\displaystyle{H_{i}=
\frac{b_{11}^{i}g_{22}-2b_{12}^{i}g_{12}+b_{22}^{i}g_{11}}{2{W}}}$ for $i=1,2$.
The Gauss curvature $K$ of $X$ in $\mathbb{E}^4_1$ is given by
\begin{equation}
\label{eq6}
K=\frac{\epsilon_{1}(b_{11}^1b_{22}^1-(b_{12}^1)^2)+
{\epsilon_{2}(b_{11}^2b_{22}^2-(b_{12}^2})^2)}{{W}},
\end{equation}
where
${\epsilon}_{1}=\langle N_{1},N_{1}\rangle$
and
${\epsilon}_{2}=\langle N_{2},N_{2}\rangle$.
When the mean curvature vector $H$ of $X$ is zero, $X$ is called as a minimal (maximal) semi-Riemann surface in $\mathbb{E}^4_1$ and when the Gaussian curvature of $X$ is zero, $X$ is called as developable (flat) semi-Riemann surface in $\mathbb{E}^{4}_{1}$. Also, $X$ is said to be a marginally trapped surface if the mean curvature vector $H$ is lightlike \cite{Chen}.

In \cite{ChenP}, the definition of the Gauss map was given as follows. Grasmanian manifold $G(2,4)$ is a space formed by all oriented 2-dimensional planes passing through the origin in $\mathbb{E}^4_{1}$. Oriented 2-dimensional planes passing through the origin in $\mathbb{E}^4_{1}$ can be defined by the unit 2-vectors. 2-vectors are elements of space $\bigwedge^2 \mathbb{E}^{4}_{1}$, that is, they are obtained with the help of wedge product $(\bigwedge)$ of vectors. The Gauss map corresponds to the oriented tangent space of semi-Riemann surface $X$ in $\mathbb{E}^{4}_{1}$ to every point of $M$. Thus, it is defined as
\begin{equation}
\nu: X \rightarrow G(2,4)\subset \mathbb{E}_{t}^{6}; \, \nu(p)=(e_{1}\wedge e_{2})(p).
\end{equation}

Now, we suppose that $X$ is a timelike surface in $\mathbb{E}^4_1$, that is, $W<0$. Thus, we can choose an orthonormal tangent frame field $e_{1},{e}_{2}$ on $M$ as below:
 \begin{equation}
 \label{eq7}
     e_{1}=\frac{1}{\sqrt{\epsilon g_{11}}}X_{u},
     \;\;\; e_{2}=\frac{1}{\sqrt{-\epsilon W g_{11}}}
     (g_{11}X_{v}-g_{12}X_{u}),
 \end{equation}
where ${\epsilon}=\langle e_{1},e_{1}\rangle=-\langle e_{2},e_{2}\rangle$. Thus, the Gauss map $\nu$ of $M$ can be given by
 \begin{equation}
 \label{eq8}
     \nu=\frac{\epsilon}{\sqrt{-W}}X_{u} \wedge X_{v}.
 \end{equation}

\section{Helicoidal Surface of Type I}

Let $\{\eta_{1},\eta_{2},\eta_{3},\eta_{4}\}$ be a standard orthonormal basis of $\mathbb{E}^4_1$, where $\eta_{1}=(1,0,0,0)$,
$\eta_{2}=(0,1,0,0)$, $\eta_{3}=(0,0,1,0)$ and $\eta_{4}=(0,0,0,1)$.
We choose as a timelike $2-$plane ${P}_1=span\{\eta_{3},\eta_{4}\}$, a hyperplane ${\Pi}_1=span\{\eta_{1},\eta_{3},\eta_{4}\}$ and a line $l_{1}=span\{\eta_{4}\}$.
Also, we suppose that 
$\beta_1:I\longrightarrow\Pi_1\subset\mathbb{E}^4_1;
\;\beta_{1}(u)=\left(x(u),0,z(u),w(u)\right)$
is a regular curve, where $x(u)\neq 0$. Thus, the parametrization of $X_1$ (called as the helicoidal surface of type I) which is obtained the rotation of the curve $\beta_1$ which leaves the timelike plane ${P}_1$ pointwise fixed followed by the translation along $l_1$ as follows: 
\begin{equation}
\label{eq9}
   X_{1}(u,v) =(x(u)\cos v,x(u)\sin v,z(u),w(u)+\lambda v),
\end{equation}
where $0\leq v<2\pi$ and $\lambda\in \mathbb{R^{+}}$. When $w$ is a constant function, $X_{1}$ is called as right helicoidal surface of type I. Also, when $z$ is a constant function, $X_{1}$ is just a helicoidal surface in $\mathbb{E}^3_{1}$.
For $\lambda=0$, the helicoidal surface which is given by \eqref{eq9} reduces to the rotational surface of elliptic type in $\mathbb{E}^4_1$ (see \cite{Dursun} and \cite{Bektas}).

By a direct calculation, we get the induced metric of $X_1$ given as follows.
\begin{equation}
\label{type1eq2}
ds^2_{X_1}=(x'^2(u)+z'^2(u)-w'^2(u))du^2-2\lambda w'(u) dudv+
(x^2(u)-\lambda^2)dv^2.
\end{equation}
with ${W}=(x^2(u)-\lambda^2)(x'^2(u)+z'^2(u))-x^2(u)w'^2(u)<0$ for all $u\in I\subset\mathbb{R}$. 
Then, we choose an orthonormal frame field 
$\{e_{1},e_{2}, N_{1}, N_{2}\}$ on $X_1$ in $\mathbb{E}^4_1$ 
such that $e_1, e_2$ are tangent to $X_1$ and $N_1, N_2$ are normal to $X_1$ as follows.
\begin{align}
\label{type1eq3}
\begin{split}
e_{1}&=\frac{1}{\sqrt{\epsilon g_{11}}}{{X_1}_u},\;\;\;\;
e_{2}=
\frac{1}{\sqrt{-\epsilon Wg_{11}}}(g_{11}{X_1}_v-{g_{12}{X_1}_u}),\\
N_{1}&=\frac{1}{\sqrt{x'^2+z'^2}}(z'\cos v, z'\sin v, -x',0),\\
N_{2}&=\frac{1}{\sqrt{-W(x'^2+z'^2)}}
(xx'w'\cos v-\lambda(x'^2+z'^2)\sin v,
xx'w'\sin v+\lambda(x'^2+z'^2)\cos v,\\
&xz'w', x(x'^2+z'^2))
\end{split}
\end{align}
where $\langle e_{1},e_{1}\rangle=-\langle e_{2},e_{2}\rangle=\epsilon=\pm 1$ 
and $\langle N_{1},N_{1}\rangle=\langle N_{2},N_{2}\rangle=1$.
For $\epsilon=1$, the surface $X_1$ has a spacelike meridian curve. 
Otherwise, it has a timelike meridian curve. 
By direct computations, we get the coefficients of the second fundamental form given as follows.
\begin{align}
\label{type1eq4}
    \begin{split}
     &b_{11}^1=\frac{x'' z'- x' z''}{\sqrt{{x'^2}+z'^2}},\;\;\; 
     b_{12}^1=b_{21}^1=0,\;\;\;
     b_{22}^1=-\frac{x z'}{\sqrt{x'^2+z'^2}},\\
     &b_{11}^2=\frac{x(w'(x'x''+z'z'')-w''(x'^2+z'^2))}
    {\sqrt{-W(x'^2+z'^2)}},\;\; 
     b_{12}^2=b_{21}^2=\frac{\lambda x'\sqrt{x'^2+z'^2}}{\sqrt{-W}},\\
    &b_{22}^2=-\frac{x^2x'w'}{\sqrt{-W(x'^2+z'^2)}}.
    \end{split}
\end{align}
Thus, the mean curvature vector $H^{X_1}$ of $X_1$ in $\mathbb{E}^4_1$ as 
\begin{equation}
H^{X_1}=H_1^{X_1} N_1 + H_2^{X_1} N_2,
\end{equation}
where $N_1, N_2$ are normal vector fields in \eqref{type1eq3},
$H_1^{X_1}$ and $H_2^{X_1}$ are given by 
\begin{align}
\label{type1eq5} 
\begin{split}
&H_{1}^{X_{1}}=\frac{(x^2-\lambda^2)(x''z'-x'z'')
-xz'(x'^2+z'^2-w'^2)}
{2W\sqrt{x'^2+z'^2}},\\
&H_{2}^{X_{1}}=\frac{x'w'(2\lambda^2-x^2)(x'^2+z'^2)+x^2x'w'^3-x(x^2-\lambda^2)(x'(x'w''-x''w')+z'(z'w''-w'z''))}
{2\sqrt{-W^3(x'^2+z'^2)}}.
\end{split}
\end{align}

\subsection{Bour's Theorem and the Gauss map for helicoidal surface of type I}
In this section, we study on Bour's theorem for timelike helicoidal surface of type I in $\mathbb{E}^4_1$ 
and we analyse the Gauss maps of isometric pair of surfaces. 

\begin{theorem}
\label{type1Bour}
A timelike helicoidal surface of type I in $\mathbb{E}^4_{1}$ given by \eqref{eq9} is isometric 
to one of the following timelike rotational surfaces in $\mathbb{E}^4_1$:
\begin{itemize}
  \item[(i)] 
\begin{equation}
\label{type1izR1}
R_{1}^1(u,v)=\left(
\begin{array}{c}
\sqrt{x^2(u)-\lambda^2}\cos{\left(v-\int{\frac{\lambda w'(u)}{x^2(u)-\lambda^2}du}\right)}\\
\sqrt{x^2(u)-\lambda^2}\sin{\left(v-\int{\frac{\lambda w'(u)}{x^2(u)-\lambda^2}du}\right)}\\
\int{\frac{a(u)x(u)x'(u)}{\sqrt{x^2(u)-\lambda^2}}du}\\
\int{\frac{b(u)x(u)x'(u)}{\sqrt{x^2(u)-\lambda^2}}du}
\end{array}
\right)
\end{equation}
so that spacelike helices on the timelike helicoidal surface of type I correspond to parallel spacelike circles 
on the timelike rotational surfaces, 
where $a(u)$ and $b(u)$ are differentiable functions satisfying the following equation:
\begin{equation}
\label{type1eq7R1}
a^2(u)-b^2(u)=
\frac{x^2(u)(z'^2(u)-w'^2(u))-\lambda^2(x'^2(u)+z'^2(u))}
{x^2(u)x'^2(u)}
\end{equation}
for all $u\in I_1$ with $x'(u)\not=0$, 

\item[(ii)] 
\begin{equation}
\label{type1izR2}
R_1^2(u,v)=\left(
\begin{array}{c}
\int{\frac{a(u)x(u)x'(u)}{\sqrt{x^2(u)-\lambda^2}}du}\\
\int{\frac{b(u)x(u)x'(u)}{\sqrt{x^2(u)-\lambda^2}}du}\\
\sqrt{x^2(u)-\lambda^2}\sinh{\left(v-\int{\frac{\lambda w'(u)}{x^2(u)-\lambda^2}du}\right)}\\
\sqrt{x^2(u)-\lambda^2}\cosh{\left(v-\int{\frac{\lambda w'(u)}{x^2(u)-\lambda^2}du}\right)}\\
\end{array}
\right)
\end{equation}
so that spacelike helices on the timelike helicoidal surface of type I correspond to parallel spacelike hyperbolas 
on the timelike rotational surfaces, 
where $a(u)$ and $b(u)$ are differentiable functions satisfying the following equation:
\begin{equation}
\label{type1eq7R2}
a^2(u)+b^2(u)=
\frac{(x^2(u)-\lambda^2)(x'^2(u)+z'^2(u))+x^2(u)(x'^2(u)-w'^2(u))}
{x^2(u)x'^2(u)}
\end{equation}
for all $u\in I_1$ with $x'(u)\not=0$,

\item[(iii)] 
\begin{equation}
\label{type1izR3}
R_1^3(u,v)=\left(
\begin{array}{c}
-\int{\frac{a(u)x(u)x'(u)}{\sqrt{\lambda^2-x^2(u)}}du}\\
-\int{\frac{b(u)x(u)x'(u)}{\sqrt{\lambda^2-x^2(u)}}du}\\
\sqrt{\lambda^2-x^2(u)}\cosh{\left(v+\int{\frac{\lambda w'(u)}{\lambda^2-x^2(u)}du}\right)}\\
\sqrt{\lambda^2-x^2(u)}\sinh{\left(v+\int{\frac{\lambda w'(u)}{\lambda^2-x^2(u)}du}\right)}\\
\end{array}
\right)
\end{equation}
so that timelike helices on the timelike helicoidal surface of type I correspond to parallel timelike hyperbolas 
on the timelike rotational surfaces,
where $a(u)$ and $b(u)$ are differentiable functions satisfying the following equation:
\begin{equation}
\label{type1eq7R3}
a^2(u)+b^2(u)=
\frac{(\lambda^2-x^2(u))(x'^2(u)+z'^2(u))-x^2(u)(x'^2(u)-w'^2(u))}
{x^2(u)x'^2(u)}
\end{equation}
for all $u\in I_2\subset\mathbb{R}$ with $x'(u)\not=0$.
\end{itemize}
\end{theorem}

\begin{proof}
Assume that $X_1$ is a timelike helicoidal surface of type I in $\mathbb{E}^4_1$ defined by \eqref{eq9}. 
Then, we have the induced metric of $X_1$ given by \eqref{type1eq2}. 
Now, we will find new coordinates $\bar{u},\bar{v}$ such that the metric becomes 
\begin{equation}
    ds^2_{X_1}=F(\bar{u})d{\bar u}^2+G(\bar{u})d\bar{v}^2,
\end{equation}
where $F({\bar{u}})$ and $G(\bar{u})$ are smooth functions. 
Set $\bar{u}=u$ and $\displaystyle{\bar{v}=v-\int{\frac{\lambda w'(u)}{x^2(u)-\lambda^2}du}}$.
Since Jacobian $\displaystyle{\frac{\partial(\bar{u},\bar{v})}{\partial(u,v)}}$ is nonzero,
it follows that $\{\bar{u}, \bar{v}\}$ are new parameters of $X_1$.  
According to the new parameters, 
the equation \eqref{type1eq2} becomes
\begin{equation}
\label{type1eq10}
 ds^2 _{X_1}
=\left(x'^2(u)+z'^2(u)-w'^2(u)
-\frac{\lambda^2w'^2(u)}{x^2(u)-\lambda^2}\right)du^2 + (x^2(u)-\lambda^2)d\overline{v}^2.
\end{equation}
Define the two subsets $I_1=\{u\in I\;|\;x^2(u)-\lambda^2>0\}$ and $I_2=\{u\in I\;|\;x^2(u)-\lambda^2<0\}$
of $I$. Then, we consider the following cases.

\textit{Case(i.)} Assume that $I_1$ is dense in the interval $I$. First,
we consider a timelike rotational surface $R_1$ in $\mathbb{E}^4_{1}$ given by 
\begin{equation}
\label{type1rtk}
    R_{1}(k,t)=(n(k)\cos{t},n(k)\sin{t},s(k),r(k))
\end{equation}
whose the induced metric is 
\begin{equation}
\label{type1eq11}
 ds^2 _{R_1}=(\dot{n}^2(k)+\dot{s}^2(k)-\dot{r}^2(k))dk^2 + n^2(k)dt^2.
\end{equation}
with $n(k)>0$. Comparing the equations \eqref{type1eq10} and \eqref{type1eq11}, 
we take $\bar{v}=t$ and $n(k)=\sqrt{x^2(u)-\lambda^2}$ and we also have
\begin{equation}
\label{type1eq12}
 \left(x'^2(u)+z'^2(u)-w'^2(u)-\frac{\lambda^2 w'^2(u)}{x^2(u)-\lambda^2}\right)du^2
=(\dot{n}^2(k)+\dot{s}^2(k)-\dot{r}^2(k)) dk^2.
\end{equation}
Set $a(u)=\frac{\dot{s}(k)}{\dot{n}(k)}$ and $b(u)=\frac{\dot{r}(k)}{\dot{n}(k)}$.
Then, we obtain 
\begin{equation}
  s=\int{\frac{a(u)x(u)x'(u)}
  {\sqrt{x^2(u)-\lambda^2}}du},\;\;\;\;
  r=\int{\frac{b(u)x(u)x'(u)}
  {\sqrt{x^2(u)-\lambda^2}}du}.
\end{equation}
Thus, we get an isometric timelike rotational surface $R_1^1$ given by \eqref{type1izR1}
satisfying \eqref{type1eq7R1}.
It can be easily seen that a spacelike helix on $X_1$ which is defined by $u=u_0$ for a constant $u_0$
corresponds to the parallel spacelike circle on $R_1^1$ 
lying on the plane $\{x_3=c_3, x_4=c_4\}$ with the radius $\sqrt{x^2_0-\lambda^2}$ for constants $c_3$ and $c_4$,
i.e., $R_{1}^1(u_0,v)=(\sqrt{x^2_0-\lambda^2}\cos{v}, \sqrt{x^2_0-\lambda^2}\sin{v},c_3,c_4)$.

Secondly, we consider a timelike rotational surface $R_{2a}$ in $\mathbb{E}^4_{1}$ given by 
\begin{equation}
\label{type2rtk}
    R_{2a}(k,t)=(n(k),p(k),r(k)\sinh{t},r(k)\cosh{t})
\end{equation}
with the induced metric as
\begin{equation}
\label{type2eq11}
 ds^2 _{R_{2a}}=(\dot{n}^2(k)+\dot{p}^2(k)-\dot{r}^2(k))dk^2 + r^2(k)dt^2
\end{equation}
with $r(k)>0$. 
Similarly, from the equations \eqref{type1eq10} and \eqref{type2eq11}, 
we take $\bar{v}=t$, $r(k)=\sqrt{x^2(u)-\lambda^2}$ 
and we have
\begin{equation}
\label{type1eq12a}
 \left(x'^2(u)+z'^2(u)-w'^2(u)-\frac{\lambda^2 w'^2(u)}{x^2(u)-\lambda^2}\right)du^2
=(\dot{n}^2(k)+\dot{p}^2(k)-\dot{r}^2(k)) dk^2.
\end{equation}
If we set $a(u)=\frac{\dot{n}(k)}{\dot{r}(k)}$ and 
$b(u)=\frac{\dot{p}(k)}{\dot{r}(k)}$, then we find
\begin{equation}
n=\int{\frac{a(u)x(u)x'(u)}
{\sqrt{x^2(u)-\lambda^2}}du},\;\;\;\;
p=\int{\frac{b(u)x(u)x'(u)}
{\sqrt{x^2(u)-\lambda^2}}du}.
\end{equation}
Thus, we get an isometric timelike rotational surface $R_1^2$ given by \eqref{type1izR2} satisfying 
\eqref{type1eq7R2}. 
It can be easily seen that a spacelike helix on $X_1$ corresponds to the parallel spacelike hyperbola 
lying on the plane $\{x_1=c_1, x_2=c_2\}$ for constants $c_1$ and $c_2$, i.e., 
$R_1^2(u_0,v)=(c_1,c_2,\sqrt{x^2_0-\lambda^2}\sinh{v},\sqrt{x^2_0-\lambda^2}\cosh{v})$.

\textit{Case(ii.)} Assume that $I_2$ is dense in the interval $I$.
We consider a timelike rotational surface $R_{2b}$ in $\mathbb{E}^4_1$ given by 
\begin{equation}
\label{type2brtk}
R_{2b}(k,t)=(n(k),p(k),s(k)\cosh{t},s(k)\sinh{t})
\end{equation}    
with the induced metric as
\begin{equation}
\label{type2eq11a}
 ds^2 _{R_{2b}}=(\dot{n}^2(k)+\dot{p}^2(k)+\dot{s}^2(k))dk^2 - s^2(k)dt^2.
\end{equation}
Considering the equations \eqref{type1eq10} and \eqref{type2eq11a}, 
we take $\bar{v}=t$, $s(k)=\sqrt{\lambda^2-x^2(u)}$ 
and we have
\begin{equation}
\label{type1eq13a}
\left(x'^2(u)+z'^2(u)-w'^2(u)-\frac{\lambda^2 w'^2(u)}{x^2(u)-\lambda^2}\right)du^2
=(\dot{n}^2(k)+\dot{p}^2(k)+\dot{s}^2(k)) dk^2.
\end{equation}
Set $a(u)=\frac{\dot{n}(k)}{\dot{s}(k)}$ and $b(u)=\frac{\dot{p}(k)}{\dot{s}(k)}$.
Then, we find
\begin{equation}
n=-\int{\frac{a(u)x(u)x'(u)}
  {\sqrt{\lambda^2-x^2(u)}}du},\;\;\;\;
p=-\int{\frac{b(u)x(u)x'(u)}
  {\sqrt{\lambda^2-x^2(u)}}du}.
\end{equation}
Thus, we get an isometric timelike rotational surface $R_1^3$ 
given by \eqref{type1izR3} satisfying \eqref{type1eq7R3}.
It can be easily seen that a timelike helix on $X_1$
corresponds to the parallel timelike hyperbola 
lying on the plane $\{x_1=c_1, x_2=c_2\}$ for constants $c_1$ and $c_2$, i.e, 
$R_1^3(u_0,v)=(c_1,c_2,\sqrt{\lambda^2-x^2_0}\cosh{v},
\sqrt{\lambda^2-x^2_0}\sinh{v})$.
\end{proof} 

Now, we find the Gauss maps of the surfaces given in Theorem \ref{type1Bour}.

\begin{lemma}
\label{Gaussmpstype1}
Let $X_1, R_1^1, R_1^2$ and $R_1^3$ be timelike surfaces in $\mathbb{E}^4_1$ given
by \eqref{eq9}, \eqref{type1izR1}, \eqref{type1izR2} and \eqref{type1izR3}, respectively. 
Then, the Gauss maps of them are given by the followings
\begin{align}
\label{type1GH}
\nu_{X_{1}}&=
\frac{\epsilon}{\sqrt{-W}}\Bigg(xx'\eta_{12}+xz'\sin{v}\eta_{13}+\left(\lambda x'\cos{v}+xw'\sin{v}\right)\eta_{14}-xz'\cos{v}\eta_{23}\notag\\
&+\left(\lambda x'\sin{v}-xw'\cos{v}\right)\eta_{24}+\lambda z'\eta_{34}\Bigg),\\
\label{type1GR11}
\nu_{R_{1}^1}&=\frac{\epsilon xx'}{\sqrt{-W}}
\Bigg(\eta_{12}+a\sin\left(v-\int\frac{\lambda w'}{x^2-\lambda^2}du\right)\eta_{13}
+b\sin\left(v-\int\frac{\lambda w'}{x^2-\lambda^2}du\right)\eta_{14}\notag\\
&-a\cos\left(v-\int\frac{\lambda w'}{x^2-\lambda^2}du\right)\eta_{23}
-b\cos\left(v-\int\frac{\lambda w'}{x^2-\lambda^2}du\right)\eta_{24}\Bigg),\\
\label{type1GR12}
\nu_{R_{1}^2}&=\frac{\epsilon xx'}{\sqrt{-W}}
\Bigg(a\cosh\left(v-\int\frac{\lambda w'}{x^2-\lambda^2}du\right)\eta_{13}
+a\sinh\left(v-\int\frac{\lambda w'}{x^2-\lambda^2}du\right)\eta_{14}\notag\\
&+b\cosh\left(v-\int\frac{\lambda w'}{x^2-\lambda^2}du\right)\eta_{23}
+b\sinh\left(v-\int\frac{\lambda w'}{x^2-\lambda^2}du\right)\eta_{24}-\eta_{34}\Bigg),\\
\label{type1GR13}
\nu_{R_{1}^3}&=-\frac{\epsilon xx'}{\sqrt{-W}}
\Bigg(a\sinh\left(v-\int\frac{\lambda  w'}{x^2-\lambda^2}du\right)\eta_{13}
+a\cosh\left(v-\int\frac{\lambda  w'}{x^2-\lambda^2}du\right)\eta_{14}\notag\\
&+b\sinh\left(v-\int\frac{\lambda  w'}{x^2-\lambda^2}du\right)\eta_{23}
+b\cosh\left(v-\int\frac{\lambda  w'}{x^2-\lambda^2}du\right)\eta_{24}+\eta_{34}\Bigg),
\end{align}
where $\{\eta_{1},\eta_{2},\eta_{3},\eta_{4}\}$ is the standard 
orthonormal bases of $\mathbb{E}_1^4$ and 
$\eta_{ij}=\eta_i\wedge \eta_j$ for $i,j=1,2,3,4$. 
\end{lemma}

\begin{proof}
Assume that $X_1$ is a timelike helicoidal surface of type I in $\mathbb{E}^4_1$ given by \eqref{eq9}.
From a direct computation, we find the Gauss map of $X_1$ by using the equation \eqref{type1eq3} in \eqref{eq8}. 
Similarly, we obtain the Gauss maps of $R_1^1, R_1^2$ and $R_1^3$ given by \eqref{type1GR11}, \eqref{type1GR12} and \eqref{type1GR13}. 
\end{proof}

For later use, we give the following lemma related to the components of the mean curvature vector
of the timelike rotational surface $R_1^1$ in $\mathbb{E}^4_1$ given by \eqref{type1izR1}. 

\begin{lemma}
Let $R_1^1$ be a timelike rotational surface in $\mathbb{E}^4_1$ defined by \eqref{type1izR1}.
Then, the mean curvature vector $H^{R_1^1}$ of $R_1^1$ in $\mathbb{E}^4_1$ 
is 
\begin{equation}
H^{R_1^1}=H_1^{R_1^1} N_1+H_2^{R_1^1} N_2, 
\end{equation}
where $N_1, N_2$ are normal vector fields in \eqref{type1eq3},
$H_1^{R_1^1}$ and $H_2^{R_1^1}$ are given by 
\begin{align}
\label{type1izm1}
H_1^{R_1^1}&=\frac{(\lambda^2-x^2)a'+axx'(b^2-a^2-1)}
{2xx'(1+a^2-b^2)\sqrt{(1+a^2)(x^2-\lambda^2)}},\\
\label{type1izm2}
H_2^{R_{1}^{1}}&=
\frac{(x^2-\lambda^2)(aa'b-b'-a^2b')-bxx'(a^2-b^2+1)}{2xx'
\sqrt{(1+a^2)(x^2-\lambda^2)(b^2-a^2-1)^{3}}}.  
\end{align}
\end{lemma}

\begin{proof}
It follows from a direct computation. 
\end{proof}

Then, we consider isometric surfaces according to Bour's theorem whose Gauss maps are same. 
\begin{theorem}
\label{type1thm2}
Let $X_{1}, R_{1}^1,R_{1}^2,R_{1}^3 $ be a timelike helicoidal surface of type I and timelike rotational surfaces in 
$\mathbb{E}^4_1$ given by \eqref{eq9}, \eqref{type1izR1}, \eqref{type1izR2} and \eqref{type1izR3}, respectively.
Then, we have the following statements.
\begin{itemize}
    \item[(i.)] If the Gauss maps of $X_1$ and $R_1^1$ are same, 
     then they are hyperplanar and minimal. Then, the parametrizations of $X_1$ and $R_1^1$ 
are given by 
\begin{equation}
\label{type1izaGH}
   X_{1}(u,v) =(x(u)\cos v,x(u)\sin v,c_1,w(u)+\lambda v)
\end{equation}
and
\begin{equation}
\label{izodönel1t}
  R_{1}^1(u,v)=\left(
\begin{array}{c}
\sqrt{x^2(u)-\lambda^2}\cos{\left(v-\int{\frac{\lambda w'(u)}{x^2(u)-\lambda^2}du}\right)}\\ \sqrt{x^2(u)-\lambda^2}\sin{\left(v-\int{\frac{\lambda w'(u)}{x^2(u)-\lambda^2}du}\right)}\\
c_2\\
\frac{1}{\sqrt{-c_3}}\arcsin{\sqrt{c_3(\lambda^2-x^2(u))}}+c_4
\end{array}
\right),
\end{equation}
where $c_1, c_2, c_3, c_4$ are arbitrary constants with $c_3<0$ and  
\begin{equation}
\label{eqcor1}
w(u)=\pm\Bigg(\sqrt{\frac{c_3\lambda^2-1}{c_3}}
\arcsin\left({\sqrt{c_3(\lambda^2-x^2(u))}}\right)-\frac{\lambda\sqrt{1-c_3\lambda^2}}{\sqrt{1+c_3\lambda^2}}\arctan{\left(
\sqrt{\frac{(1+c_3\lambda^2)(x^2(u)-\lambda^2)}{\lambda^2(1-c_3(x^2(u)-\lambda^2))}}\right)}\Bigg).
\end{equation}

\item[(ii.)] The Gauss maps of $X_1$ and $R_1^2$ or $R_1^3$ are definitely different. 
\end{itemize}
\end{theorem}

\begin{proof}
Assume that $X_{1}$ is a timelike helicoidal surface of type I in $\mathbb{E}^4_{1}$ defined by \eqref{eq9} and $R_{1}^1,R_{1}^2,R_{1}^3$ are timelike rotational surfaces in $\mathbb{E}^4_1$ defined by \eqref{type1izR1}, \eqref{type1izR2} and \eqref{type1izR3}, respectively.  
From Lemma \ref{Gaussmpstype1}, 
we know the Gauss maps of $X_1, R_1^1, R_1^2$ and $R_1^3$ given by \eqref{type1GH}, \eqref{type1GR11}, \eqref{type1GR12} 
and \eqref{type1GR13}, respectively.
Then, we consider the Gauss maps of each surfaces.\\ 
\textit{(i)}
Suppose that $X_1$ and $R_1^1$ have the same Gauss maps.
From \eqref{type1GH} and \eqref{type1GR11}, we get the following system of equations:
\begin{align}
  \label{type1eq12t}
    xz'\sin{v}&=a
    xx'\sin{\left(v-\int\frac{\lambda w'}{x^2-\lambda^2}du\right)},\\
  \label{type1eq13}
    xz'\cos{v}&=a
    xx'\cos{\left(v-\int\frac{\lambda w'}{x^2-\lambda^2}du\right)},\\
  \label{type1eq14}
    \lambda x'\cos{v}+xw'\sin{v}&=bxx'\sin{\left(v-\int
    \frac{\lambda w'}{x^2-\lambda^2}du\right)},\\
  \label{type1eq15}
    \lambda x'\sin{v}-xw'\cos{v}&=-bxx'\cos{\left(v-\int
    \frac{\lambda w'}{x^2-\lambda^2}du\right)},\\
  \label{type1eq16}
    \lambda z'&=0.
\end{align}
Due to $\lambda \neq 0$, the equation \eqref{type1eq16} gives 
$z'(u)=0$. Then, from
the equations \eqref{type1eq12t} and \eqref{type1eq13}
we get $a(u)=0$. 
Therefore, it can be easily seen that 
the timelike surfaces $X_{1}$ and $R_{1}^1$ are hyperplanar,
that is, they are lying in $\mathbb{E}^3_1$.
Moreover, the equations \eqref{type1eq5} and \eqref{type1izm1} 
imply that $H_{1}^{X_{1}}=H_{1}^{R_{1}^1}=0$.
Also, from the equations \eqref{type1eq5} and \eqref{type1izm2}, we have
\begin{align}
\label{type1eq18}
\begin{split}
H_{2}^{X_{1}}&=\frac{x'^2w'(2\lambda^2-x^2)
+x^2w'^3+x(x^2-\lambda^2)(x''w'-x'w'')}{2(x^2w'^2-x'^2(x^2-\lambda^2))^{3/2}},\\
H_{2}^{R_{1}^1}&=
\frac{bxx'(b^2-1)-b'(x^2-\lambda^2)}{2xx'\sqrt{(x^2-\lambda^2)(b^2-1)^{3}}}.
\end{split}
\end{align}
Using $z'(u)=a(u)=0$, from the equation \eqref{type1eq7R1} we have 
\begin{equation}
 \label{type1eq18a}   
 b^2=\frac{x^2w'^2+\lambda^2 x'^2}{x^2x'^2}. 
\end{equation}
Using the equation \eqref{type1eq18a} in \eqref{type1eq18}, 
we get
\begin{equation}
 \label{type1eq18b}
H_{2}^{R_{1}^1}=\frac{x^2w'(x'^2w'(2\lambda^2-x^2)+x^2w'^3
+x(x^2-\lambda^2)(x''w'-x'w''))}{2(x^2w'^2-x'^2(x^2-\lambda^2))^{3/2}\sqrt{(x^2w'^2+
\lambda^2x'^2)(\lambda^2-x^2)}}.
\end{equation}
Thus, we get
$\displaystyle{H_2^{R_1^1}
=\frac{x^2w'}
{\sqrt{(x^2w'^2+\lambda^2x'^2)(x^2-\lambda^2)}}H_2^{X_1}}$. 
Moreover, using equations \eqref{type1eq14} and \eqref{type1eq15}, we obtain the following equations
\begin{align}
\label{type1eq21}
    xw'=bxx'\cos{\left(\int\frac{\lambda w'}{x^2-\lambda^2}du\right)},\\
\label{type1eq22}
    \lambda x'=-bxx'\sin{\left(\int\frac{\lambda w'}{x^2-\lambda^2}du\right)}.
\end{align}
Considering the equations \eqref{type1eq21} and \eqref{type1eq22} together, we have
\begin{equation}
\label{type1eq23}
    \frac{xw'}
    {\lambda x'}=-\cot{\left(\int\frac{\lambda w'}{x^2-\lambda^2}du\right)}.
\end{equation}
Taking the derivative of \eqref{type1eq23} with respect to $u$, 
we find 
\begin{equation}
\label{type1eq24}
    \lambda^2(xx'w''+w'(2x'^2-
    xx''))+x^2(w'(w'^2 -x'^2)+x(x''w'-x'w''))=0
\end{equation}
which implies $H_{2}^{X_{1}}=H_{2}^{R_{1}^1}=0$. 
Thus, we get the desired results.
Since $R_1$ is minimal, from the equation \eqref{type1izm1} 
we have the following differential equation
\begin{equation}
\label{coreq2}
 (x^2-\lambda^2)b'+xx'b=xx'b^3  
\end{equation}
which is a Bernoulli equation. 
Then, the general solution of this equation is found as 
\begin{equation}
\label{coreq3}
   b^2=\frac{1}{1+c_3(x^2-\lambda^2)}   
\end{equation}
for an arbitrary negative constant $c_3$.
Comparing the equations \eqref{type1eq18a} and \eqref{coreq3}, we get
\begin{equation}
 w(u)
 =\pm\sqrt{1-c_3\lambda^2}
 \int{\frac{x'(u)}{x(u)}\sqrt{\frac{x^2(u)-\lambda^2}{1+c_3(x^2(u)-\lambda^2)}}}du 
\end{equation}
whose solution is given by \eqref{eqcor1} for $c_3<0$. 
Moreover, using the last component of $R_1(u,v)$ in \eqref{izodönel1t}, 
we have
\begin{equation}
\int\frac{x(u)x'(u)}{\sqrt{(x^2(u)-\lambda^2)(1+c_3(x^2(u)-\lambda^2))}}du
=\pm\frac{1}{\sqrt{-c_3}}\arcsin{\sqrt{-c_3(x^2(u)-\lambda^2)}}+c_4
\end{equation}
for any arbitrary constant $c_4$.
\\
\textit{(ii.)}
Suppose that $X_1$ and $R_1^1$ have the same Gauss maps. 
Comparing the equations \eqref{type1GH} and \eqref{type1GR12}, we get $x(u)=0$ or $x'(u)=0$ which give 
$\nu_{R_{1}^2}=0$. That is a contradiction. 
Thus, their Gauss maps are definitely different.
Similarly, we show that the Gauss maps of the $X_{1}$ and $R_{1}^{3}$ surfaces are definitely different.
\end{proof}

\begin{remark}
Ikawa studied Bour's theorem for helicoidal surfaces in $\mathbb{E}^3_1$ and 
he also established the parametrizations of the isometric surfaces when they have the same Gauss map.
Taking $x(u)=u$ in Theorem \ref{type1thm2}, we get the cases obtained in \cite{Ikawa2}. 
Moreover, he determined the minimal rotational surfaces in $\mathbb{E}^3_1$, \cite{Ikawa2}. 
The rotational surface given by \eqref{izodönel1t} has the same form of surface in Proposition 3.4, \cite{Ikawa2}.
\end{remark}

\begin{remark}
If $w'(u)=0$ for $u\in I\subset\mathbb{R}$, 
then the timelike helicoidal surface given by \eqref{eq9} reduces to the timelike right helicoidal surface in $\mathbb{E}^4_1$. 
On the other hand, $W=(x^2(u)-\lambda^2)(x'^2(u)+z'^2(u))<0$ for $x^2(u)-\lambda^2<0$. Thus, 
from Theorem \ref{type1Bour}, we get the timelike rotational surfaces $R_1^3(u,v)$ 
which are isometric to the timelike right helicoidal surface in $\mathbb{E}^4_1$. 
Also, Theorem \ref{type1thm2} implies that the Gauss maps of such surfaces are definitely different. 
\end{remark}

Now, we give an example by using Theorem \ref{type1thm2}. 
 
\begin{example}
If we choose $x(u)=u$, $\lambda=1$, $c_3 =-1/2$ and $c_4=0$, then isometric surfaces in \eqref{type1izaGH} and \eqref{izodönel1t} are given as follows
\begin{equation*}
 X_{1}(u,v)=\left(u\cos v, u\sin v,\sqrt{3}\left(\arcsin{\sqrt{\frac{u^2-1}{2}}}
-\arctan{\sqrt{\frac{u^2-1}{u^2+1}}}\right)+v\right)   
\end{equation*}
and
\begin{equation*}
\label{type1izaGR}
R_{1}^1(u,v)=\left(
\begin{array}{c}
\sqrt{u^2-1}\cos{\left(v-\frac{1}{2}
\arctan{\left(\frac{2u^2-3}{\sqrt{3}\sqrt{-u^4+4u^2-3}}\right)}\right)}\\ 
\sqrt{u^2-1}\sin{\left(v-\frac{1}{2}
\arctan{\left(\frac{2u^2-3}{\sqrt{3}\sqrt{-u^4+4u^2-3}}\right)}\right)}\\
\sqrt{2}\arcsin{\sqrt{\frac{u^2-1}{2}}}
\end{array}
\right).    
\end{equation*}
For $1.32\leq u\leq\ 1.72$ and $0\leq v<2\pi$,  
the graphs of timelike helicoidal surface $X_{1}$ and timelike rotational surface $R_{1}$ in $\mathbb{E}^3_1$ 
can be plotted by using Mathematica 10.4 as follows:

\begin{figure}[h]
\centering
\begin{subfigure}{.3\textwidth}
  \centering
  \includegraphics[width=5cm, height=3.5cm]{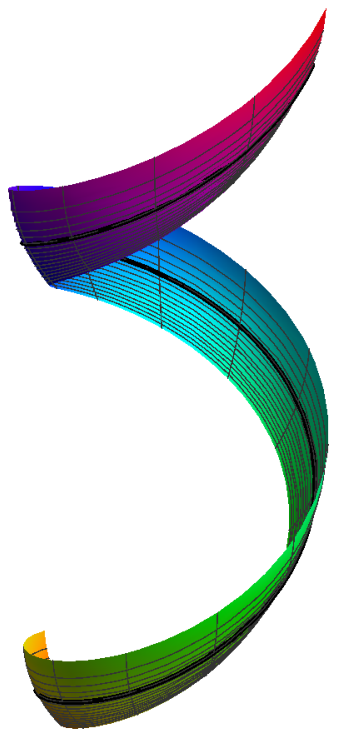}
  \caption{}
  \label{fig:sub1}
\end{subfigure}%
\begin{subfigure}{.3\textwidth}
  \centering
  \includegraphics[width=5cm, height=3.5cm]{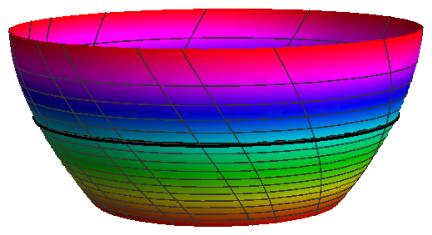}
  \caption{}
  \label{fig:sub2}
\end{subfigure}
\caption{(A) Timelike helicoidal surface of type I; spacelike helix and (B)  Timelike rotational surface; spacelike circle.}
\end{figure}
\end{example}

\section{Helicoidal Surface of Type IIa}
Let us choose a timelike $2-$plane $P_{2}=span\{\eta_{1},\eta_{2}\}$, a hyperplane $\Pi_{2a}=span\{\eta_{1},\eta_{2},\eta_{4}\}$ and a line $l_{2}=span\{\eta_{1}\}$. Also, we suppose that 
$\beta_{2a}:I\longrightarrow\Pi_{2a}\subset\mathbb{E}^4_1;
\;\beta_{2a}(u)=\left(x(u),y(u),0, w(u)\right)$
is a regular curve, where $w(u)\neq 0.$ Thus, the parametrization of $X_{2a}$ (called as the helicoidal surface of type IIa) which is obtained the rotation of the curve $\beta_{2a}$ which leaves 
the timelike plane ${P}_{2}$ pointwise fixed followed by the translation along $l_{2}$ as follows:
\begin{equation}
\label{eq10}
X_{2a}(u,v)=(x(u)+\lambda v,y(u),w(u)\sinh v,w(u)\cosh v),
\end{equation}
where, $v\in\mathbb{R}$ and $\lambda\in \mathbb{R^{+}}$. When $x$ is a constant function, $X_{2a}$ is called as right helicoidal surface of type IIa. Also, when $y$ is a constant function, $X_{2a}$ is just a helicoidal surface in $\mathbb{E}^3_{1}$.
For $\lambda=0$, the helicoidal surface which is given by \eqref{eq10} reduces to the rotational surface of hyperbolic type in $\mathbb{E}^4_1$ (see \cite{Dursun} and \cite{Bektas}).

By a direct calculation, we get the induced metric of $X_{2a}$ given as follows
\begin{equation}
\label{type2mtrc1}
 ds^2_{X_{2a}}=(x'^2(u)+y'^2(u)-w'^2(u))du^2+2\lambda x'(u)dudv + (\lambda^2+w^2(u))dv^2
\end{equation}
with ${W}=(\lambda^2+w^2(u))(y'^2(u)-w'^2(u))+x'^2(u)w^2(u)<0$.
Then, we choose an orthonormal frame field 
$\{e_{1},e_{2}, N_{1}, N_{2}\}$ on $X_{2a}$ in $\mathbb{E}^4_1$ 
such that $e_1, e_2$ are tangent to $X_{2a}$ and $N_1, N_2$ are normal to $X_{2a}$ as follows.
\begin{align}
\label{type2eq3}
\begin{split}
e_{1}&=\frac{1}{\sqrt{\epsilon g_{11}}}X_{{2a}_u},\;\;\;\;
e_{2}=\frac{1}{\sqrt{-\epsilon Wg_{11}}}(g_{11}X_{{2a}_v}-{g_{12}X_{{2a}_u}}),\\
N_{1}&=\frac{1}{\sqrt{w'^2-y'^2}}(0,w', y'\sinh v, y'\cosh v),\\
N_{2}=&\frac{1}{\sqrt{-W(w'^2-y'^2)}}
(-w(w'^2-y'^2),-wx'y',
\lambda(w'^2-y'^2)\cosh v-x'ww'\sinh v,\\
&\lambda(w'^2-y'^2)\sinh v-x'ww'\cosh v)
\end{split}
\end{align}
where $\langle e_{1},e_{1}\rangle=-\langle e_{2},e_{2}\rangle=\epsilon$ 
and $\langle N_{1},N_{1}\rangle=\langle N_{2},N_{2}\rangle=1$. 
For $\epsilon=1$, the surface $X_{2a}$ has a spacelike meridian curve. 
Otherwise, it has a timelike meridian curve. 
By direct computations, we get the coefficients of the second fundamental form given as follows.
\begin{align}
\label{type2eq4}
    \begin{split}
     &b_{11}^1=\frac{y'' w'-y' w''}{\sqrt{w'^2-y'^2}},\;\;\; 
     b_{12}^1=b_{21}^1=0,\;\;\;
     b_{22}^1=-\frac{w y'}{\sqrt{w'^2-y'^2}},\\
     &b_{11}^2=\frac{w(x'(w'w''-y'y'')+x''(y'^2-w'^2))}
    {\sqrt{W(y'^2-w'^2)}},\;\; 
     b_{12}^2=b_{21}^2=\frac{\lambda w'\sqrt{w'^2-y'^2}}{\sqrt{-W}},\\
     &b_{22}^2=\frac{x'w^2w'}{\sqrt{W(y'^2-w'^2)}}.
    \end{split}
\end{align}
Thus, the mean curvature vector $H^{X_{2a}}$ of $X_{2a}$ in $\mathbb{E}^4_1$ as 
\begin{equation}
H^{X_{2a}}=H_1^{X_{2a}} N_1 + H_2^{X_{2a}} N_2,
\end{equation}
where $N_1, N_2$ are normal vector fields in \eqref{type2eq3},
$H_1^{X_{2a}}$ and $H_2^{X_{2a}}$ are given by 

\begin{align}
\label{type2eq5} 
\begin{split}
&H_{1}^{X_{2a}}=\frac{(\lambda^2+w^2)(y''w'-y'w'')-w y'(x'^2+y'^2-w'^2)}{2W\sqrt{w'^2-y'^2}},\\
&H_{2}^{X_{2a}}=\frac{x'w'(2\lambda^2+w^2)(y'^2-w'^2)+ x'^3w^2w'+w(\lambda^2+w^2)(x''(y'^2-w'^2)+x'(w'w''-y'y''))}
{2\sqrt{W^3(y'^2-w'^2)}}.
\end{split}
\end{align}

\subsection{Bour's Theorem and the Gauss map for helicoidal surfaces IIa}
In this section, we study on Bour's theorem for timelike helicoidal surface of type IIa in $\mathbb{E}^4_1$ 
and we analyse the Gauss maps of isometric pair of surfaces.  

\begin{theorem}
\label{type2thm1}
A timelike helicoidal surface of type IIa in $\mathbb{E}^4_{1}$ given by \eqref{eq10} is isometric to one of the following timelike rotational surfaces in $\mathbb{E}^4_1$
\begin{itemize}
\item[(i)]
 \begin{equation}
\label{type2izR1}
  R_{2a}^1(u,v)=\left(
\begin{array}{c}
\sqrt{\lambda^2+w^2(u)}\cos{\left(v+\int{\frac{\lambda x'(u)}{\lambda^2+w^2(u)}du}\right)}\\ \sqrt{\lambda^2+w^2(u)}\sin{\left(v+\int{\frac{\lambda x'(u)}{\lambda^2+w^2(u)}du}\right)}\\\int{\frac{a(u)w(u)w'(u)}{\sqrt{\lambda^2+w^2(u)}}du}\\
\int{\frac{b(u)w(u)w'(u)}{\sqrt{\lambda^2+w^2(u)}}du}
\end{array}
\right)
\end{equation}
so that spacelike helices on the timelike helicoidal surface of type IIa correspond to parallel spacelike circles on the timelike rotational surface, where $a(u)$ and $b(u)$ are differentiable functions satisfying the following equation:
\begin{equation}
\label{type2aeq7R1}
a^2(u)-b^2(u)=
\frac{\lambda^2(y'^2(u)-w'^2(u))+w^2(u)(x'^2(u)+y'^2(u)-2w'^2(u))}
{w^2(u)w'^2(u)}
\end{equation}

\item[(ii)]
\begin{equation}
\label{type2izR2}
  R_{2a}^2(u,v)=\left(
\begin{array}{c}
\int{\frac{a(u)w(u)w'(u)}{\sqrt{\lambda^2+w^2(u)}}du}\\
\int{\frac{b(u)w(u)w'(u)}{\sqrt{\lambda^2+w^2(u)}}du}\\
\sqrt{\lambda^2+w^2(u)}\sinh{\left(v+\int{\frac{\lambda x'(u)}{\lambda^2+w^2(u)}du}\right)}\\ \sqrt{\lambda^2+w^2(u)}\cosh{\left(v+\int{\frac{\lambda x'(u)}{\lambda^2+w^2(u)}du}\right)}
\end{array}
\right)
\end{equation}
so that spacelike helices on the timelike helicoidal surface of type IIa correspond to parallel spacelike hyperbolas on the timelike rotational surface, where $a(u)$ and $b(u)$ are differentiable functions satisfying the following equation: 
\begin{equation}
 \label{type2eq7a}
    a^2(u)+b^2(u)
    =\frac{w^2(u)(x'^2(u)+y'^2(u))+\lambda^2(y'^2(u)-w'^2(u))}{w^2(u)w'^2(u)}.
 \end{equation}
\end{itemize}
\end{theorem}

\begin{proof}
Assume that $X_{2a}$ is a timelike helicoidal surface of type IIa in $\mathbb{E}^4_1$ defined by \eqref{eq10}. 
Then, we have the induced metric of $X_{2a}$ given by \eqref{type2mtrc1}. 
Now, we will find new coordinates $\bar{u},\bar{v}$ such that the metric becomes 
\begin{equation}
    ds^2_{X_{2a}}=F(\bar{u})du^2+G(\bar{u})d\bar{v}^2,
\end{equation}
where $F({\bar{u}})$ and $G(\bar{u})$ are smooth functions. 
Set $\bar{u}=u$ and $\displaystyle{\bar{v}=v+\int{\frac{\lambda x'(u)}{\lambda^2+w^2(u)}du}}$.
Since Jacobian $\displaystyle{\frac{\partial(\bar{u},\bar{v})}{\partial(u,v)}}$ is nonzero,
it follows that $\{\bar{u}, \bar{v}\}$ are new parameters of $X_1$. 
According to the new parameters, 
the equation \eqref{type2mtrc1} becomes
\begin{equation}
\label{type2eq10}
ds^2 _{X_{2a}}=\left(x'^2(u)+y'^2(u)-w'^2(u)-\frac{\lambda^2 x'^2(u)}{\lambda^2+w^2(u)}\right)du^2 
+ (\lambda^2+w^2(u))d\overline{v}^2.
\end{equation}
Then, we consider the following cases.

First, we consider a timelike rotational surface $R_1$ in $\mathbb{E}^4_1$ given by \eqref{type1rtk}.
Then, we have the induced metric of $R_1$ given by \eqref{type1eq11}. 
Comparing the equations \eqref{type1eq11} and \eqref{type2eq10}, 
we take $\bar{v}=t$ and $n(k)=\sqrt{\lambda^2+w^2(u)}$ and we also have
\begin{equation}
\label{type2eq12b}
\left(x'^2(u)+y'^2(u)-w'^2(u)-\frac{\lambda^2 x'^2(u)}{\lambda^2+w^2(u)}\right)du^2
=(\dot{n}^2(k)+\dot{s}^2(k)-\dot{r}^2(k))dk^2. 
\end{equation}
Set $a(u)=\frac{\dot{s}(k)}{\dot{n}(k)}$ and 
$b(u)=\frac{\dot{r}(k)}{\dot{n}(k)}$.
Then, we obtain 
\begin{equation}
  s=\int{\frac{a(u)w(u)w'(u)}
  {\sqrt{\lambda^2+w^2(u)}}du},\;\; r=\int{\frac{b(u)w(u)w'(u)}
  {\sqrt{\lambda^2+w^2(u)}}du}.
\end{equation}
Thus, we get an isometric timelike rotational surface $R_{2a}^1$ given by \eqref{type2izR1}
satisfying \eqref{type2aeq7R1}.
It can be easily seen that a spacelike helix on $X_{2a}$ corresponds to parallel spacelike circle 
lying on the plane $\{x_3=c_3, x_4=c_4\}$ with the radius $\sqrt{\lambda^2+w^2_0}$ for constants 
$c_3$ and $c_4$, i.e., 
$R_{2a}^1(u_0,v)=(\sqrt{\lambda^2+w^2_0}\cos{v},\sqrt{\lambda^2+w^2_0}\sin{v},c_3,c_4)$.

Secondly, we consider a timelike rotational surface $R_{2a}$ in $\mathbb{E}^4_{1}$ given by
\eqref{type2rtk}. Then, we know the induced metric given by \eqref{type2eq11}.  
Comparing the equations \eqref{type2eq11} and \eqref{type2eq10}, 
we take $\bar{v}=t$ and $r(k)=\sqrt{\lambda^2+w^2(u)}$ 
and we also have
\begin{equation}
\label{type2eq12a}
 \left(x'^2(u)+y'^2(u)-w'^2(u)-\frac{\lambda^2 x'^2(u)}{\lambda^2+w^2(u)}\right)du^2=(\dot{n}^2(k)+\dot{p}^2(k)-\dot{r}^2(k))dk^2. 
\end{equation}
Set $a(u)=\frac{\dot{n}(k)}{\dot{r}(k)}$ and 
$b(u)=\frac{\dot{p}(k)}{\dot{r}(k)}$. 
Then, we obtain
\begin{equation}
  n=\int{\frac{a(u)w(u)w'(u)}
  {\sqrt{\lambda^2+w^2(u)}}du},\;\;\;\; p=\int{\frac{b(u)w(u)w'(u)}
  {\sqrt{\lambda^2+w^2(u)}}du}.
\end{equation}
Thus, we get an isometric timelike rotational surface $R_{2a}^2$ given by \eqref{type2izR2}.
It can be easily seen that a spacelike helix on $X_{2a}$ which is defined by $u=u_0$
for a constant $u_0$ corresponds to the parallel spacelike hyperbola
lying on the plane $\{x_1=c_1, x_2=c_2\}$ for constants $c_1$ and $c_2$, i.e., 
$R_{2a}^2(u_0,v)=(c_1,c_2,\sqrt{\lambda^2+w^2_0}\sinh{v},\sqrt{\lambda^2+w^2_0}\cosh{v})$.
\end{proof} 

\begin{lemma}
\label{Gaussmpstype2a}
Let $X_{2a}, R_{2a}^1$ and $R_{2a}^{2}$ be timelike surfaces in $\mathbb{E}^4_1$ given by
\eqref{eq10}, \eqref{type2izR1} and \eqref{type2izR2}, respectively.
Then, the Gauss maps of them are given by the followings
\begin{align}
\label{type2GH}
\nu_{X_{2a}}&=\frac{\epsilon}{\sqrt{-W}}\Bigg(-\lambda y'\eta_{12}+(x'w\cosh{v}-\lambda w'\sinh{v})\eta_{13}
+( x'w\sinh{v}-\lambda w'\cosh{v})\eta_{14}\notag\\
&+y'w\cosh{v}\eta_{23}   
+y'w\sinh{v}\eta_{24}-ww'\eta_{34}\Bigg),\\
\label{type2GR}
\nu_{R_{2a}^{1}}&=\frac{\epsilon ww'}{\sqrt{-W}}
\Bigg(\eta_{12}+a\sin\left(v+\int\frac{\lambda x'}{\lambda^2+w^2}du\right)\eta_{13}
+b\sin\left(v+\int\frac{\lambda x'}{\lambda^2+w^2}du\right)\eta_{14}\notag\\
& -a\cos\left(v+\int\frac{\lambda x'}{\lambda^2+w^2}du\right)\eta_{23}
-b\cos\left(v+\int\frac{\lambda x'(u)}{\lambda^2+w^2(u)}du\right)\eta_{24}\Bigg),\\
\label{type2GRa}
\nu_{R_{2a}^{2}}&=\frac{\epsilon ww'}{\sqrt{-W}}
\Bigg(a\cosh\left(v+\int\frac{\lambda x'}{\lambda^2+w^2}du\right)\eta_{13}
+a\sinh\left(v+\int\frac{\lambda x'}{\lambda^2+w^2}du\right)\eta_{14}\notag\\
&
+b\cosh\left(v+\int\frac{\lambda x'}{\lambda^2+w^2}du\right)\eta_{23} +b\sinh\left(v+\int\frac{\lambda x'}{\lambda^2+w^2}du\right)\eta_{24}-\eta_{34}\Bigg),
\end{align}
where $\{\eta_{1},\eta_{2},\eta_{3},\eta_{4}\}$ is the standard 
orthonormal bases of $\mathbb{E}_1^4$ and 
$\eta_{ij}=\eta_i\wedge \eta_j$ for $i,j=1,2,3,4$. 
\end{lemma}

\begin{proof}
Assume that $X_{2a}$ is a timelike helicoidal surface of type IIa in $\mathbb{E}^4_1$
given by \eqref{eq10}.
From a direct computation, we find the Gauss map of $X_{2a}$ by using the equation \eqref{type2eq3} in \eqref{eq8}.
Similarly, we obtain the Gauss maps of $R_{2a}^1$ and $R_{2a}^2$ given by \eqref{type2GR} and \eqref{type2GRa}. 
\end{proof}

For later use, we give the following lemma related to the components of the mean curvature vector
of the timelike rotational surface $R_{2a}^2$ given by \eqref{type2izR2}.

\begin{lemma}
Let $R_{2a}^2$ be a timelike rotational surface in $\mathbb{E}^4_1$ given by \eqref{type2izR2}. 
Then, the mean curvature vector $H^{R_{2a}^2}$ of $R_{2a}^2$ in $\mathbb{E}^4_1$ is 
\begin{equation}
H^{R_{2a}^2}=H_1^{R_{2a}^2} N_1+H_2^{R_{2a}^2} N_2,
\end{equation}
where $N_1, N_2$ are normal vector fields in \eqref{type2eq3},  
$H_1^{R_{2a}^2}$ and $H_2^{R_{2a}^2}$ are given by 
\begin{align}
\label{type2izm1}
\begin{split}
H_1^{R_{2a}^2}&=\frac{(w^2+\lambda^2)b'-(a^2+b^2-1)bww'}
{2ww'(a^2+b^2-1)\sqrt{(1-b^2)(w^2+\lambda^2)}},\\
H_2^{R_{2a}^2}&=
\frac{(w^2+\lambda^2)(a'(1-b^2)+abb')-aww'(a^2+b^2-1)}{2ww'\sqrt{(b^2-1)(w^2+\lambda^2)(a^2+b^2-1)^3}}.
\end{split}
\end{align}
\end{lemma}

\begin{proof}
It follows from a direct computation.
\end{proof}

Then, we consider isometric surfaces according to Bour's theorem whose Gauss maps are same. 

\begin{theorem}
\label{type2thm4}
Let $X_{2a}, R_{2a}^{1}$ and $R_{2a}^{2}$ be a timelike helicoidal surface of type IIa 
and timelike rotational surfaces in  $\mathbb{E}^4_{1}$ 
given by \eqref{eq10}, \eqref{type2izR1} and \eqref{type2izR2}, respectively. 
Then, we have the following statements.
\begin{itemize}
\item[(i.)] The Gauss maps of $X_{2a}$ and $R_{2a}^1$ are definitely different.

\item[(ii.)] If the surfaces $X_{2a}$ and $R_{2a}^2$ have the same Gauss maps, then they are hyperplanar and minimal. 
Then, the parametrizations of $X_{2a}$ and $R_{2a}^2$ can be explicitly determined by 
\begin{equation}
\label{type2izaGH}
   X_{2a}(u,v) =\left(x(u)+\lambda v,c_1,w(u)\sinh v,w(u)\cosh v\right)
\end{equation}
and
\begin{equation}
\label{type2izaGR}
  R_{2a}^{2}(u,v)=\left(
\begin{array}{c}
\pm\frac{1}{\sqrt{c_3}}\arcsinh{\sqrt{c_3(\lambda^2+w^2(u))}}+c_4\\
c_2\\
\sqrt{\lambda^2+w^2(u)}\sinh{\left(v+\int{\frac{\lambda x'(u)}{\lambda^2+w^2(u)}du}\right)}\\ 
\sqrt{\lambda^2+w^2(u)}\cosh{\left(v+\int{\frac{\lambda x'(u)}{\lambda^2+w^2(u)}du}\right)}
\end{array}
\right),
\end{equation}
where $c_1, c_2, c_3, c_4$ are arbitrary constants with $c_3>0$ and 
\begin{equation}
\label{eqcor2}
\begin{split}
x(u)=&\pm\left(\sqrt{1+c_3\lambda^2}
\arcsinh{\sqrt{c_3(\lambda^2+w^2(u))}}-\sqrt{c_3\lambda^2}
\arctanh{\left(\frac{\lambda \sqrt{1+c_3(\lambda^2+w^2(u))}}
{\sqrt{(1+c_3\lambda^2)(\lambda^2+w^2(u))}}\right)}\right).    
\end{split}
\end{equation}
\end{itemize}
\end{theorem}

\begin{proof}
Assume that $X_{2a}$ is a timelike helicoidal surface of type I in $\mathbb{E}^4_{1}$ 
given by \eqref{eq10} and $R_{2a}^1, R_{2a}^2$ are timelike rotational surfaces $\mathbb{E}^4_1$ 
given by \eqref{type2izR1} and \eqref{type2izR2}, respectively.
From Lemma \ref{Gaussmpstype2a}, 
we have the Gauss maps of $X_{2a}, R_{2a}^1$ and $R_{2a}^2$   given by \eqref{type2GH}, \eqref{type2GR}
and \eqref{type2GRa}, respectively.\\ 
\textit{(i.)}
Suppose that the Gauss maps of $X_{2a}$ and $R_{2a}^{1}$ are same. Then, from the equations \eqref{type2GH} and \eqref{type2GR}, 
we get $w(u)=0$ or $w'(u)=0$ which implies $\nu_{R_{2a}^1}=0$. That is a contradiction. 
Thus, their Gauss maps are definitely different.\\
\textit{(ii)}
Suppose that the surfaces $X_{2a}$ and $R_{2a}^{2}$ have the same Gauss maps. 
From \eqref{type2GH} and \eqref{type2GRa}, we get the following system of equations:
\begin{align}
\label{type2eq12}
    \lambda y'&=0,\\
\label{type2eq13}
    x'w\cosh{v}-\lambda w'\sinh{v}&=a
    ww'\cosh\left(v+\int\frac{\lambda x'}{\lambda^2+w^2}du\right),\\
\label{type2eq14}
    x'w\sinh{v}-\lambda w'\cosh{v}&=a
    ww'\sinh\left(v+\int\frac{\lambda x'}{\lambda^2+w^2}du\right),\\
\label{type2eq15}
    y'w\cosh{v}&=bww'\cosh\left(v+\int\frac{\lambda x'}{\lambda^2+w^2}du\right),\\
\label{type2eq16}
     y'w\sinh{v}&=bww'\sinh\left(v+\int\frac{\lambda x'}{\lambda^2+w^2}du\right).
\end{align}
Due to $\lambda \neq 0 $,  the equation \eqref{type2eq12} gives $y'(u)=0$.
Then, from the equations \eqref{type2eq15} and \eqref{type2eq16}
imply $b(u)=0$.
Therefore, it can be easily seen that 
the surfaces $X_{2a}$ and $R_{2a}^{2}$ are hyperplanar,
that is, they are lying in $\mathbb{E}^3_1$. 
Moreover, the equations \eqref{type2eq5} and \eqref{type2izm1} imply that $H_{1}^{X_{2a}}=H_{1}^{R_{2a}^2}=0$
and 
\begin{align}
\label{type2eq18}
\begin{split}
H_{2}^{X_{2a}}&=-\frac{x'w'^2(2\lambda^2+w^2)
-w^2x'^3+w(\lambda^2+w^2)(x''w'-x'w'')}{2(w'^2(\lambda^2+w^2)-w^2x'^2)^{3/2}},\\
H_{2}^{R_{2a}^2}&=
\frac{a'(w^2+\lambda^2)+aww'(1-a^2)}
{2ww'\sqrt{(w^2+\lambda^2)(1-a^2)^3}}.
\end{split}
\end{align}
Using $b(u)=0$, from the equation \eqref{type2eq7a} we have 
\begin{equation}
 \label{type2eq18a}   
 a^2(u)=\frac{w^2(u)x'^2(u)-\lambda^2 w'^2(u)}{w^2(u)w'^2(u)}.
\end{equation}
Using the equation \eqref{type2eq18a} in \eqref{type2eq18}, we get
\begin{equation}
 \label{type2eq18b}
H_{2}^{R_{2a}^{2}}=\frac{w^2x'(x'w'^2(2\lambda^2+w^2)
-w^2x'^3+w(\lambda^2+w^2)(x''w'-x'w''))}{2(w'^2(\lambda^2+w^2)-w^2x'^2)^{3/2}\sqrt{\left(\lambda^2+w^2\right)\left(w'^2\lambda^2-w^2x'^2\right)}}.    
\end{equation}
Which implies $H_2^{R_{2a}^{2}}
=-\frac{w^2x'}{\sqrt{(w^2x'^2-\lambda^2w'^2)(\lambda^2+w^2)}}H_2^{X_{2a}}$. 
Moreover, using equations \eqref{type2eq13} and \eqref{type2eq14}, we obtain the following equations
\begin{align}
 \label{type2eq21}
    x'w=aww'\cosh{\left(\int\frac{\lambda x'}{\lambda^2+w^2}du\right)},\\
\label{type2eq22}
    \lambda w'=-aww'\sinh{\left(\int\frac{\lambda x'}{\lambda^2+w^2}du\right)}.   
\end{align}
Considering the equations \eqref{type2eq21} and \eqref{type2eq22} together, we have
\begin{equation}
\label{type2eq23}
    -\frac{x'w}
    {\lambda w'}=\coth{\left(\int\frac{\lambda x'}{\lambda^2+w^2}du\right)}.
\end{equation}
If we take the derivative of the equation \eqref{type2eq23} with respect to $u$, \eqref{type2eq23} becomes
\begin{equation}
\label{type2eq24}
    \lambda(x'w'^2(2\lambda^2+w^2)
-w^2x'^3+w(\lambda^2+w^2)(x''w'-x'w''))=0
\end{equation}
which implies $H_{2}^{X_{2a}}=H_{2}^{R_{2a}^{2}}=0$ in the equation 
\eqref{type2eq18}. 
Now, we determine the parametrizations of the surfaces $X_{2a}$ and $R_{2a}^{2}$.
Since $R_{2a}^{2}$ is minimal, from the equation \eqref{type2eq18} 
we have the following Bernoulli differential equation
\begin{equation}
\label{coreqtype2}
 (\lambda^2+w^2)a'+ww'a=ww'a^3  
\end{equation}
whose solution is given by
\begin{equation}
\label{cortype2eq1}
   a^2=\frac{1}{1+c_3(\lambda^2+w^2)}   
\end{equation}
for an arbitrary positive constant $c_3$. 
Comparing the equations \eqref{type2eq18a} and \eqref{cortype2eq1}, we get
\begin{equation}
 x(u)
=\pm\sqrt{1+c_3\lambda^2}
 \int{\frac{w'(u)}{w(u)}
 \sqrt{\frac{w^2(u)+\lambda^2}{1+c_3(w^2(u)+\lambda^2)}}}du.  
\end{equation}
whose solution is given by \eqref{eqcor2} for $c_3>0$. 
Moreover, using the last component of $R_{2a}^2(u,v)$ in \eqref{type2izaGR}, we have 
\begin{equation}
\int\frac{w(u)w'(u)}{\sqrt{(\lambda^2+w^2(u))(1+c_3(\lambda^2+w^2(u))}}du=\pm\frac{1}{\sqrt{c_3}}\arcsinh{\sqrt{c_3(\lambda^2+w^2(u)}}+c_4
\end{equation}
for any arbitrary constant $c_4$.
\end{proof}

\begin{remark}
Ikawa studied Bour's theorem for helicoidal surfaces in $\mathbb{E}^3_1$ and 
he also established the parametrizations of the isometric surfaces when they have the same Gauss map.
Taking $w(u)=u$ in Theorem \ref{type2thm4}, we get the cases obtained in \cite{Ikawa2}. 
Moreover, he determined the minimal rotational surfaces in $\mathbb{E}^3_1$, \cite{Ikawa2}. 
The rotational surface given by \eqref{type2izaGR} has the same form of surface in Proposition 3.2, \cite{Ikawa2}.
\end{remark}

\begin{remark}
If $x'(u)=0$ for $u\in I\subset\mathbb{R}$, 
then the timelike helicoidal surface given by \eqref{eq10} reduces to the right timelike helicoidal surface in $\mathbb{E}^4_1$. 
Thus, 
from Theorem \ref{type2thm1}, 
we get the timelike rotational surfaces $R_{2a}^{1}(u,v)$ and $R_{2a}^2(u,v)$
which are isometric to the timelike right helicoidal surface in $\mathbb{E}^4_1$. 
Also, Theorem \ref{type2thm4} implies that the Gauss maps of $X_{2a}$ and $R_{2a}^{1}$ are definitely different. 
If the timelike right helicoidal surface and $R_{2a}^2$ have the same Gauss map, then we get $a^2(u)=-\frac{\lambda^2}{w^2(u)}$
which gives a contradiction. Thus, they have the different Gauss maps. 
\end{remark}

Now, we give an example by using Theorem \ref{type2thm4}. 
\begin{example}
If we choose $w(u)=u$, $\lambda=c_3 =1$ and $c_4=0$, then isometric surfaces in \eqref{type2izaGH} and \eqref{type2izaGR} are given as follows 
\begin{equation*}
\label{type2izaGHn}
X_{2a}(u,v) =\left(\sqrt{2}\arcsinh{\sqrt{1+u^2}}-\arctanh{\sqrt{\frac{2+u^2}{2+2u^2}}}+ v,u\sinh v,u\cosh v\right)
\end{equation*}
and 
\begin{equation*}
\label{type2izaGRn}
R_{2}^{2a}(u,v)=\left(
\begin{array}{c}
\arcsinh{\sqrt{1+u^2}}\\
\sqrt{1+u^2}\sinh\left(v+\ln{\frac{u}{\sqrt{3u^2+2\sqrt{2u^4+6u^2+4}+4}}} \right)\\ \sqrt{1+u^2}\cosh\left(v+\ln{\frac{u}{\sqrt{3u^2+2\sqrt{2u^4+6u^2+4}+4}}} \right)
\end{array}
\right).
\end{equation*}
For $1.19\leq u\leq\ 10$ and $-1.5\leq v<1.5$,  
the graphs of timelike helicoidal surface $X_{2a}$ and timelike rotational surface $R_{2}^2$ in $\mathbb{E}^3_1$ can be plotted by using Mathematica 10.4 as follows:

\begin{figure}[h]
\centering
\begin{subfigure}{.4\textwidth}
  \centering
  \includegraphics[width=5cm, height=4cm]{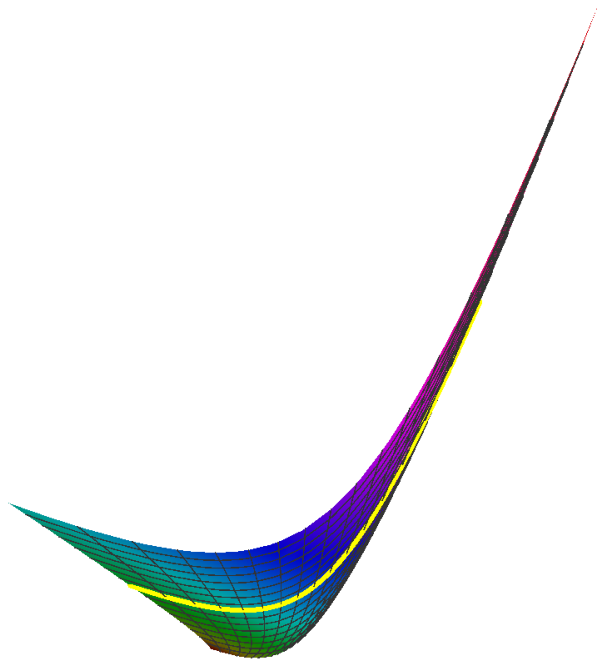}
  \caption{}
  \label{fig:sub1}
\end{subfigure}%
\begin{subfigure}{.4\textwidth}
  \centering
  \includegraphics[width=5cm, height=4cm]{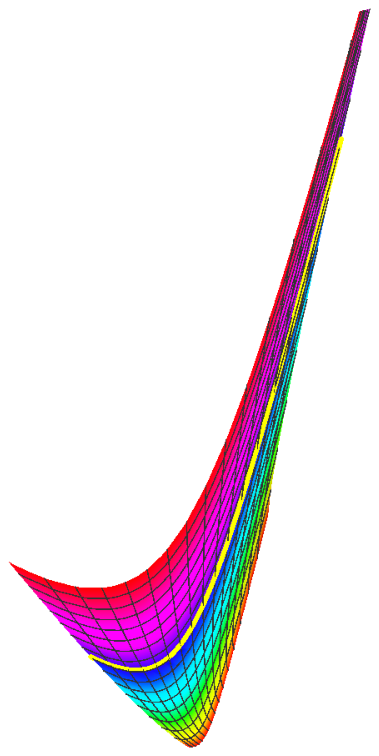}
  \caption{}
  \label{fig:sub2}
\end{subfigure}
\caption{(A) Timelike helicoidal surface of type IIa; spacelike helix and (B)  Timelike rotational surface; spacelike hyperbola.}
\end{figure}
\end{example}

\section{Helicoidal Surface of Type IIb}
Let us choose a timelike $2-$plane $P_{2}=span\{\eta_{1},\eta_{2}\}$, a hyperplane $\Pi_{2b}=\{\eta_{1},\eta_{2},\eta_{3}\}$ and a line $l_{2}=span\{\eta_{1}\}$. Also, we suppose that 
$\beta_{2b}:I\longrightarrow\Pi_{2b}\subset\mathbb{E}^4_1;
\;\beta_{2b}(u)=\left(x(u),y(u),z(u),0\right)$
is a regular curve, where $z(u)\neq 0.$ Thus, the parametrization of $X_{2b}$ (called as the timelike helicoidal surface of type IIb) which is obtained the rotation of the curve $\beta_{2b}$ which leaves 
the timelike plane ${P}_{2}$ pointwise fixed followed by the translation along $l_{2}$ as follows:
\begin{equation}
\label{eq11}
X_{2b}(u,v)=(x(u)+\lambda v,y(u),z(u)\cosh v,z(u)\sinh v),
\end{equation}
where, $v\in\mathbb{R}$ and $\lambda\in \mathbb{R^{+}}$. When $x$ is a constant function, $X_{2b}$ is called as timelike right helicoidal surface of type IIb. Also, when $y$ is a constant function, $X_{2b}$ is just a timelike helicoidal surface in $\mathbb{E}^3_{1}$.
For $\lambda=0$, the helicoidal surface which is given by \eqref{eq11} reduces to the rotational surface of hyperbolic type in $\mathbb{E}^4_1$.

By a direct calculation, we get the induced metric of $X_{2b}$ given as follows.
\begin{equation}
\label{type2bmtrc1}
ds^2 _{X_{2b}}=(x'^2(u)+y'^2(u)+z'^2(u))du^2+2\lambda x'(u)dudv + (\lambda^2-z^2(u))dv^2
\end{equation}
with ${W}=(\lambda^2-z^2(u))(y'^2(u)+z'^2(u))-x'^2(u)z^2(u)<0$.
Then, we choose an orthonormal frame field 
$\{e_{1},e_{2}, N_{1}, N_{2}\}$ on $X_{2b}$ in $\mathbb{E}^4_1$ 
such that $e_1, e_2$ are tangent to $X_{2b}$ and $N_1, N_2$ are normal to $X_{2b}$ as follows.

\begin{align}
\label{type2beq3}
\begin{split}
e_{1}&=\frac{1}{\sqrt{\epsilon g_{11}}}X_{{2b}_u},\;\;\;\;
e_{2}=\frac{1}{\sqrt{-\epsilon Wg_{11}}}(g_{11}X_{{2b}_v}-{g_{12}X_{{2b}_u}}),\\
N_{1}&=\frac{1}{\sqrt{y'^2+z'^2}}(0,-z', y'\cosh v, y'\sinh v),\\
N_{2}=&\frac{1}{\sqrt{-W(y'^2+z'^2)}}
(-z(y'^2+z'^2),zx'y',
x'zz'\cosh v-\lambda(y'^2+z'^2)\sinh v,\\
&x'zz'\sinh v-\lambda(y'^2+z'^2)\cosh v)
\end{split}
\end{align}
where $\langle e_{1},e_{1}\rangle=-\langle e_{2},e_{2}\rangle=\epsilon$ 
and $\langle N_{1},N_{1}\rangle=\langle N_{2},N_{2}\rangle=1$. 
For $\varepsilon=1$, the surface $X_{2b}$ has a spacelike meridian curve.
Otherwise, it has a timelike meridian curve.

By direct computations, we get the coefficients of the second fundamental form given as follows.
\begin{align}
\label{type2beq4}
    \begin{split}
     &b_{11}^1=\frac{y' z''-y'' z'}{\sqrt{y'^2+z'^2}},\;\;\; 
     b_{12}^1=b_{21}^1=0,\;\;\;
     b_{22}^1=\frac{z y'}{\sqrt{y'^2+z'^2}},\\
     &b_{11}^2=\frac{z(x'(y'y''+z'z'')-x''(y'^2+z'^2))}
    {\sqrt{-W(y'^2+z'^2)}},\;\; 
     b_{12}^2=b_{21}^2=\frac{\lambda z'\sqrt{y'^2+z'^2}}{\sqrt{-W}},\\
     &b_{22}^2=\frac{x'z^2z'}{\sqrt{-W(y'^2+z'^2)}}.
    \end{split}
\end{align}
Thus, the mean curvature vector $H^{X_{2b}}$ of $X_{2b}$ in $\mathbb{E}^4_1$ as 
\begin{equation}
H^{X_{2b}}=H_1^{X_{2b}} N_1 + H_2^{X_{2b}} N_2,
\end{equation}
where $N_1, N_2$ are normal vector fields in \eqref{type2beq3},
$H_1^{X_{2b}}$ and $H_2^{X_{2b}}$ are given by 
\begin{align}
\label{type2beq5} 
\begin{split}
H_{1}^{X_{2b}}&=\frac{(\lambda^2-z^2)(y'z''-z'y'')+zy'(x'^2+y'^2+z'^2)}
{2W\sqrt{y'^2+z'^2}},\\
H_{2}^{X_{2b}}&=\frac{1}
{2\sqrt{-W^3(y'^2+z'^2)}}\Bigg(x'z'((z^2-2\lambda^2)(y'^2+z'^2)+z^2(x'^2-zz''))\\
&+\lambda^2z(x'(z'z''+y'y'')-x''(y'^2+z'^2))+z^3(x''(y'^2+z'^2)-x'y'y'')\Bigg).
\end{split}
\end{align}

\subsection{Bour's Theorem and the Gauss map for helicoidal surfaces of type IIb}
In this section, we study on Bour's theorem for timelike helicoidal surface of type IIb 
in $\mathbb{E}^4_1$ and we analyse the Gauss maps of isometric pair of surfaces. 

\begin{theorem}
\label{type2bBour}
A timelike helicoidal surface of type IIb in $\mathbb{E}^4_{1}$ given by \eqref{eq11} is isometric to one of 
the following timelike rotational surfaces in $\mathbb{E}^4_1$:
\begin{itemize}
\item[(i)] 
\begin{equation}
\label{type2bizR1}
R_{2b}^1(u,v)=\left(
\begin{array}{c}
\sqrt{\lambda^2-z^2(u)}\cos{\left(v+\int{\frac{\lambda x'(u)}{\lambda^2-z^2(u)}du}\right)}\\
\sqrt{\lambda^2-z^2(u)}\sin{\left(v+\int{\frac{\lambda x'(u)}{\lambda^2-z^2(u)}du}\right)}\\
-\int{\frac{a(u)z(u)z'(u)}{\sqrt{\lambda^2-z^2(u)}}du}\\
-\int{\frac{b(u)z(u)z'(u)}{\sqrt{\lambda^2-z^2(u)}}du}
\end{array}
\right)
\end{equation}
so that spacelike helices on the timelike helicoidal surface of type IIb correspond to parallel spacelike circles 
on the timelike rotational surfaces,
where $a(u)$ and $b(u)$ are differentiable functions satisfying the following equation:
\begin{equation}
\label{type2eqR1a}
a^2(u)-b^2(u)=
\frac{\lambda^2(y'^2(u)+z'^2(u))-z^2(u)(x'^2(u)+y'^2(u)+2z'^2(u))}
{z^2(u)z'^2(u)}
\end{equation}
for all $u\in I_1\subset\mathbb{R}$.

\item[(ii)] 
\begin{equation}
\label{type2bizR2}
R_{2b}^2(u,v)=\left(
\begin{array}{c}
-\int{\frac{a(u)z(u)z'(u)}{\sqrt{\lambda^2-z^2(u)}}du}\\
-\int{\frac{b(u)z(u)z'(u)}{\sqrt{\lambda^2-z^2(u)}}du}\\\
\sqrt{\lambda^2-z^2(u)}\sinh{\left(v+\int{\frac{\lambda x'(u)}{\lambda^2-z^2(u)}du}\right)}\\
\sqrt{\lambda^2-z^2(u)}\cosh{\left(v+\int{\frac{\lambda x'(u)}{\lambda^2-z^2(u)}du}\right)}
\end{array}
\right)
\end{equation}
so that spacelike helices on the timelike helicoidal surface of type IIb correspond to parallel spacelike hyperbolas
on the timelike rotational surfaces,
where $a(u)$ and $b(u)$ are differentiable functions satisfying the following equation:
\begin{equation}
\label{type2eqR2a}
a^2(u)+b^2(u)=
\frac{\lambda^2(y'^2(u)+z'^2(u))-z^2(u)(x'^2(u)+y'^2(u))}
{z^2(u)z'^2(u)}
\end{equation}
for all $u\in I_1\subset\mathbb{R}$,

\item[(iii)] 
\begin{equation}
\label{type2bizR3}
R_{2b}^3(u,v)=\left(
\begin{array}{c}
\int{\frac{a(u)z(u)z'(u)}{\sqrt{z^2(u)-\lambda^2}}du}\\
\int{\frac{b(u)z(u)z'(u)}{\sqrt{z^2(u)-\lambda^2}}du}\\\sqrt{z^2(u)-\lambda^2}\cosh{\left(v-\int{\frac{\lambda x'(u)}{z^2(u)-\lambda^2}du}\right)}\\
\sqrt{z^2(u)-\lambda^2}\sinh{\left(v-\int{\frac{\lambda x'(u)}{z^2(u)-\lambda^2}du}\right)}
\end{array}
\right)
\end{equation}
so that timelike helices on the timelike helicoidal surface of type IIb correspond to parallel timelike hyperbolas 
on the timelike rotational surfaces,
where $a(u)$ and $b(u)$ are differentiable functions satisfying the following equation:
\begin{equation}
\label{type2eqR3a}
a^2(u)+b^2(u)=
\frac{z^2(u)(x'^2(u)+y'^2(u))-\lambda^2(y'^2(u)+z'^2(u))}
{z^2(u)z'^2(u)}
\end{equation}
with $z'(u)\not=0$ for all $u\in I_2\subset\mathbb{R}$.
\end{itemize}
\end{theorem}

\begin{proof}
Assume that $X_{2b}$ is a timelike helicoidal surface of type IIb in $\mathbb{E}^4_1$ defined by \eqref{eq11}. 
Then, we have the induced metric of $X_{2b}$ given by \eqref{type2bmtrc1}. 
Now, we will find new coordinates $\bar{u},\bar{v}$ such that the metric becomes 
\begin{equation}
    ds^2_{X_{2b}}=F(\bar{u})du^2+G(\bar{u})d\bar{v}^2,
\end{equation}
where $F({\bar{u}})$ and $G(\bar{u})$ are smooth functions. 
Set $\bar{u}=u$ and $\overline{v}=v+\int{\frac{\lambda x'(u)}{\lambda^2-z^2(u)}du}$.
Since Jacobian $\displaystyle{\frac{\partial(\bar{u},\bar{v})}{\partial(u,v)}}$ is nonzero,
it follows that $\{\bar{u}, \bar{v}\}$ are new parameters of $X_1$. 
According to the new parameters, 
the equation \eqref{type2bmtrc1} becomes
\begin{equation}
\label{type2beq10}
 ds^2 _{X_{2b}}=\left(x'^2(u)+y'^2(u)+z'^2(u)+\frac{\lambda^2 x'^2(u)}{z^2(u)-\lambda^2}\right)du^2 
+ (\lambda^2-z^2(u))d\overline{v}^2.
\end{equation}
Define two subsets $I_1=\{u\in I\;|\; z^2(u)-\lambda^2<0\}$ and $I_2=\{u\in I\;|\; z^2(u)-\lambda^2>0\}$.
Then, we consider the following cases.\\
\textit{Case(i)} Assume that $I_1$ is dense in $I$. 
First, we consider a timelike rotational surface $R_1$ in $\mathbb{E}^4_1$ given by \eqref{type1rtk}. 
Comparing the equations \eqref{type1eq11} and \eqref{type2beq10}, 
we take $\bar{v}=t$ and $n(k)=\sqrt{\lambda^2-z^2(u)}$
and we also have
\begin{equation}
\label{type2beq12b}
\left(x'^2(u)+y'^2(u)+z'^2(u)+\frac{\lambda^2 x'^2(u)}{z^2(u)-\lambda^2}\right)du^2
=(\dot{n}^2(k)+\dot{s}^2(k)-\dot{r}^2(k))dk^2.
\end{equation}
Set $a(u)=\frac{\dot{s}(k)}{\dot{n}(k)}$ and $b(u)=\frac{\dot{r}(k)}{\dot{n}(k)}$. 
Then, we obtain
\begin{equation}
  s=-\int{\frac{a(u)z(u)z'(u)}
  {\sqrt{\lambda^2-z^2(u)}}du},\;\;\;\; r=-\int{\frac{b(u)z(u)z'(u)}
  {\sqrt{\lambda^2-z^2(u)}}du}.
\end{equation}
Thus, we get an isometric timelike rotational surface $R_{2b}^1$ given by \eqref{type2bizR1}
satisfying \eqref{type2eqR1a}.
It can be easily seen that a spacelike helix on $X_{2b}$ which is defined by $u=u_0$
for a constant $u_0$ corresponds to the parallel spacelike circle
lying on the plane $\{x_3=c_3, x_4=c_4\}$ with the radius $\sqrt{\lambda^2-z^2_0}$ for constants 
$c_3$ and $c_4$, i.e., 
$R_{2b}^1(u_0,v)=(\sqrt{\lambda^2-z_0^2}\cos{v},\sqrt{\lambda^2-z_0^2}\sin{v},c_3,c_4)$.

Secondly, we consider a timelike rotational surface $R_{2a}$ in $\mathbb{E}^4_{1}$ 
given by \eqref{type2rtk}.
Then, we have the equation \eqref{type2eq11}. 
Comparing the equations \eqref{type2eq11} and \eqref{type2beq10}, we take
$\bar{v}=t$ and $r(k)=\sqrt{\lambda^2-z^2(u)}$ and we also have 
\begin{align}
\label{type2beq1Ga}
\left(x'^2(u)+y'^2(u)+z'^2(u)+\frac{\lambda^2 x'^2(u)}{z^2(u)-\lambda^2}\right)du^2
=(\dot{n}^2(k)+\dot{p}^2(k)-\dot{r}^2(k))dk^2. 
\end{align}
Set $a(u)=\frac{\dot{n}(k)}{\dot{r}(k)}$ and $b(u)=\frac{\dot{p}(k)}{\dot{r}(k)}$. 
We find
\begin{equation}
n=-\int{\frac{a(u)z(u)z'(u)}
{\sqrt{\lambda^2-z^2(u)}}du},\;\;\;\; 
p=-\int{\frac{b(u)z(u)z'(u)}{\sqrt{\lambda^2-z^2(u)}}du}.
\end{equation}
Thus, we get an isometric timelike rotational surface $R_{2b}^2$ given by \eqref{type2bizR2} 
satisfying \eqref{type2eqR2a}. 
It can be easily seen that a spacelike helix on $X_{2b}$ which is defined by $u=u_0$
for a constant $u_0$ corresponds to
parallel spacelike hyperbolas lying on the plane $\{x_1=c_1, x_2=c_2\}$ for constants $c_1$ and $c_2$
i.e., $R_{2b}^2(u_0,v)=(c_1,c_2,\sqrt{\lambda^2-z_0^2}\sinh{v},\sqrt{\lambda^2-z_0^2}\cosh{v})$.

\textit{Case (ii)} Assume that $I_2$ is dense in $I$. 
Then, we consider a timelike rotational surface $R_{2b}$ in $\mathbb{E}^4_{1}$ given by \eqref{type2brtk}. 
Comparing the equations \eqref{type2eq11a} and \eqref{type2beq10}, 
we take $\bar{v}=t$ and $s(k)=\sqrt{z^2(u)-\lambda^2}$ and we also have
\begin{equation}
\label{type2beq12a}
\left(x'^2(u)+y'^2(u)+z'^2(u)+\frac{\lambda^2 x'^2(u)}{z^2(u)-\lambda^2}\right)du^2
=(\dot{n}^2(k)+\dot{p}^2(k)+\dot{s}^2(k))dk^2. 
\end{equation}
Set $a(u)=\frac{\dot{n}(k)}{\dot{s}(k)}$ and $b(u)=\frac{\dot{p}(k)}{\dot{s}(k)}$. 
Then, we obtain
\begin{equation}
  n=\int{\frac{a(u)z(u)z'(u)}
  {\sqrt{z^2(u)-\lambda^2}}du},\;\; p=\int{\frac{b(u)z(u)z'(u)}
  {\sqrt{z^2(u)-\lambda^2}}du}.
\end{equation}
Thus, we get the timelike isometric rotational surface $R_{2b}^3$ given by \eqref{type2bizR3}
satisfying \eqref{type2eqR3a}.
It can be easily seen that a timelike helix on $X_{2b}$ corresponds 
to the parallel timelike hyperbolas lying on the plane $\{x_1=c_1, x_2=c_2\}$ for constants $c_1$ and $c_2$, 
i.e., $R_{2b}^3(u_0,v)=(c_1,c_2,\sqrt{z_0^2-\lambda^2}\cosh{v},\sqrt{z_0^2-\lambda^2}\sinh{v})$.
\end{proof}

\begin{lemma}
\label{Gaussmpstype2b}
Let $X_{2b}, {R_{2b}^{1}}, R_{2b}^{2}$ and ${R_{2b}^{3}}$ be timelike surfaces in $\mathbb{E}^4_1$ given by 
\eqref{eq11}, \eqref{type2bizR1}, \eqref{type2bizR2} and \eqref{type2bizR3}, respectively.
The Gauss maps of them are given by
\begin{align}
\label{type2bGH}
\nu_{X_{2b}}&=
\frac{\epsilon}{\sqrt{-W}}\Bigg(-\lambda y'\eta_{12}+\left(x'z\sinh{v}-\lambda z'\cosh{v}\right)\eta_{13}
+\left(x'z\cosh{v}-\lambda z'\sinh{v}\right)\eta_{14}\notag
\\&+y'z\sinh{v}\eta_{23} +y'z\cosh{v}\eta_{24}+zz'\eta_{34}\Bigg),\\
\label{type2bGRa}
\nu_{R_{2b}^{1}}&=-\frac{\epsilon zz'}{\sqrt{-W}}
\Bigg(\eta_{12}+a\sin\left(v+\int\frac{\lambda x'}{\lambda^2-z^2}du\right)\eta_{13}
+b\sin\left(v+\int\frac{\lambda x'}{\lambda^2-z^2}du\right)\eta_{14}\notag\\
& -a\cos\left(v+\int\frac{\lambda x'}{\lambda^2-z^2}du\right)\eta_{23}
-b\cos\left(v+\int\frac{\lambda x'}{\lambda^2-z^2}du\right)\eta_{24}\Bigg),\\
\label{type2bGRb}
\nu_{R_{2b}^{2}}&=-\frac{\epsilon zz'}{\sqrt{-W}}
\Bigg(a\cosh\left(v+\int\frac{\lambda x'}{\lambda^2-z^2}du\right)\eta_{13}
+a\sinh\left(v+\int\frac{\lambda x'}{\lambda^2-z^2}du\right)\eta_{14}\notag\\
&+b\cosh\left(v+\int\frac{\lambda x'}{\lambda^2-z^2}du\right)\eta_{23}
+b\sinh\left(v+\int\frac{\lambda x'}{\lambda^2-z^2}du\right)\eta_{24}-\eta_{34}\Bigg),
\end{align}
\begin{align}
\label{type2bGRc}
\nu_{R_{2b}^{3}}&=\frac{\epsilon zz'}{\sqrt{-W}}
\Bigg(a\sinh\left(v-\int\frac{\lambda x'}{z^2-\lambda^2}du\right)\eta_{13}
+a\cosh\left(v-\int\frac{\lambda x'}{z^2-\lambda^2}du\right)\eta_{14}\notag\\
&+b\sinh\left(v-\int\frac{\lambda x'}{z^2-\lambda^2}du\right)\eta_{23}
+b\cosh\left(v-\int\frac{\lambda x'}{z^2-\lambda^2}du\right)\eta_{24}+\eta_{34}\Bigg),
\end{align}
where $\{\eta_{1},\eta_{2},\eta_{3},\eta_{4}\}$ is the standard 
orthonormal bases of $\mathbb{E}_1^4$ and 
$\eta_{ij}=\eta_i\wedge \eta_j$ for $i,j=1,2,3,4$. 
\end{lemma}

\begin{proof}
Assume that $X_{2b}$ is a timelike helicoidal surface of type IIb in $\mathbb{E}^4_1$
given by \eqref{eq11}.
From a direct computation, we find the Gauss map of $X_{2b}$ by using the equation \eqref{type2beq3} in \eqref{eq8}.
Similarly, we obtain the Gauss maps of ${R_{2b}^{1}}, R_{2b}^{2}$ and ${R_{2b}^{3}}$ given by \eqref{type2GR} and \eqref{type2GRa}. 
\end{proof}

For later use, we find the components of the mean curvature vector
of the timelike rotational surface $R_{2b}^3$ given by \eqref{type2bizR3} as follows.

\begin{lemma}
Let $R_{2b}^2$ and $R_{2b}^3$ be timelike rotational surfaces in $\mathbb{E}^4_1$ 
given by \eqref{type2bizR2} and \eqref{type2bizR3}.
\begin{itemize}
\item[(i.)] The mean curvature vector $H^{R_{2b}^2}$ of $R_{2b}^2$ in $\mathbb{E}^4_1$ is
\begin{equation}
 H^{R_{2b}^2}=H_1^{R_{2b}^2} N_1+H_2^{R_{2b}^2} N_2, 
\end{equation}
where $N_1, N_2$ are normal vector fields in \eqref{type2beq3},
$H_1^{R_{2b}^2}$ and $H_2^{R_{2b}^2}$ are given by 
\begin{align}
\label{type2bizm1}
\begin{split}
H_1^{R_{2b}^2}&=\frac{b'(z^2-\lambda^2)-bzz'(a^2+b^2-1)}
{2zz'(a^2+b^2-1)\sqrt{(b^2-1)(z^2-\lambda^2)}},\\
H_2^{R_{2b}^2}&=
\frac{(z^2-\lambda^2)(a'(b^2-1)-abb')+azz'(a^2+b^2-1)}{2zz'\sqrt{(b^2-1)(\lambda^2-z^2)(a^2+b^2-1)^3}}.
\end{split}
\end{align}

\item[(ii.)] The mean curvature vector $H^{R_{2b}^3}$ of $R_{2b}^3$ in $\mathbb{E}^4_1$ is 
\begin{equation}
 H^{R_{2b}^3}=H_1^{R_{2b}^3} N_1+H_2^{R_{2b}^3} N_2, 
\end{equation}
where $N_1, N_2$ are normal vector fields in \eqref{type2beq3},
$H_1^{R_{2b}^3}$ and $H_2^{R_{2b}^3}$ are given by
\begin{align}
\label{type2bizm2}
\begin{split}
H_1^{R_{2b}^3}&=\frac{b'(z^2-\lambda^2)+bzz'(a^2+b^2+1)}
{2zz'\sqrt{(b^2+1)(z^2-\lambda^2)(a^2+b^2+1)}},\\
H_2^{R_{2b}^3}&=
\frac{(z^2-\lambda^2)(a'(1+b^2)-abb')+azz'(a^2+b^2+1)}{2zz'\sqrt{(b^2+1)(z^2-\lambda^2)(a^2+b^2+1)^3}}.
\end{split}
\end{align}
\end{itemize}
\end{lemma}

\begin{proof}
It follows from a direct calculation. 
\end{proof}

Then, we consider isometric surfaces according to Bour's theorem whose Gauss maps are same. 

\begin{theorem}
\label{type2bthrm6}
Let $X_{2b}, R_{2b}^1, R_{2b}^2$ and $R_{2b}^3$ be a timelike helicoidal surface of type IIb and 
timelike rotational surfaces in $\mathbb{E}^4_1$ given by \eqref{eq11},
\eqref{type2bizR1}, \eqref{type2bizR2} and \eqref{type2bizR3}, respectively. 
Then, we have the following statements.
\begin{itemize}
\item[(i.)] The Gauss maps of $X_{2b}$ and $R_{2b}^1$ are definitely different. 

\item[(ii.)]If the Gauss maps of the surfaces $X_{2b}$ and $R_{2b}^2$ are same, 
then they are hyperplanar and minimal. 
Then, the parametrizations of the surfaces $X_{2b}$ and $R_{2b}^2$ can be explicitly determined by
\begin{equation}
\label{type2bizaGH}
   X_{2b}(u,v) =(x(u)+\lambda v,c_1,z(u)\cosh v,z(u)\sinh v)
\end{equation}
and
\begin{equation}
\label{type2bizaGR}
R_{2b}^{2}(u,v)=\left(
\begin{array}{c}
\pm\frac{1}{\sqrt{c_3}}\arcsinh{\sqrt{c_3 (\lambda^2-z^2(u))}+c_4}\\
c_2\\
\sqrt{\lambda^2-z^2(u)}\sinh{\left(v+\int{\frac{\lambda x'(u)}{\lambda^2-z^2(u)}du}\right)}\\ 
\sqrt{\lambda^2-z^2(u)}\cosh{\left(v+\int{\frac{\lambda x'(u)}{\lambda^2-z^2(u)}du}\right)}\\
\end{array}
\right),
\end{equation}
where $c_1, c_2, c_3$ and $c_4$ are arbitrary constants with $c_3>0$ and 
\begin{equation}
\label{eqcor2b}
\begin{split}
x(u)&=\pm\Bigg(\frac{\sqrt{1+c_3\lambda^2}}{\sqrt{c_3}}
\arcsinh{\sqrt{c_3(\lambda^2-z^2(u))}}\\
&+\lambda\arctanh\left(\frac{\sqrt{(1+c_3\lambda^2)(\lambda^2-z^2(u))}}{\lambda \sqrt{1+c_3\lambda^2}}\right)\Bigg),    
\end{split}
\end{equation}
where $z(u)\not=0$ and $z'(u)\not=0$.

\item[(iii)] If the Gauss maps of the surfaces $X_{2b}$ and $R_{2b}^3$ are same, 
then they are hyperplanar and minimal. 
Then, the parametrizations of the surfaces $X_{2b}$ and $R_{2b}^3$ can be explicitly determined by
\begin{equation}
\label{type2bizaGH}
   X_{2b}(u,v) =\left(x(u)+\lambda v,c_1,z(u)\cosh v,z(u)\sinh v\right)
\end{equation}
and
\begin{equation}
\label{type2bizaGR}
  R_{2b}^{3}(u,v)=\left(
\begin{array}{c}
\pm\frac{1}{\sqrt{c_3}}\arccosh{\sqrt{c_3\left(z^2(u)-\lambda^2
\right)}+c_4}\\
c_2\\
\sqrt{z^2(u)-\lambda^2}\cosh{\left(v-\int{\frac{\lambda x'(u)}{z^2(u)-\lambda^2}du}\right)}\\ \sqrt{z^2(u)-\lambda^2}\sinh{\left(v-\int{\frac{\lambda x'(u)}{z^2(u)-\lambda^2}du}\right)}\\
\end{array}
\right),
\end{equation}
where $c_1, c_2, c_3$ and $c_4$ are arbitrary constants with $c_3>0$ and 
\begin{align}
\label{eqcor2c}
\begin{split}
x(u)&=\pm\frac{1}{\sqrt{2}}\Bigg(\sqrt{1+c_3\lambda^2}
\arcsinh{\sqrt{c_3(z^2(u)-\lambda^2)-1}}\\
&-\lambda\sqrt{c_3}
\arctanh{\left(\frac{\lambda \sqrt{c_3(z^2(u)-\lambda^2)-1}}
{\sqrt{(1+c_3\lambda^2)(z^2(u)-\lambda^2)}}\right)}\Bigg).    
\end{split}
\end{align}
\end{itemize}
\end{theorem}

\begin{proof}
Assume that $X_{2b}$ is a timelike helicoidal surface of type I in $\mathbb{E}^4_{1}$ given by \eqref{eq11} and $R_{2}^1, R_{2}^2, R_{2}^3$ are timelike rotational surfaces $\mathbb{E}^4_1$ given by \eqref{type2bizR1}, \eqref{type2bizR2} and \eqref{type2bizR3}, respectively.
From Lemma \ref{Gaussmpstype2b}, 
we have the Gauss maps of $X_{2b}, R_{2b}^1, R_{2b}^2$ and $R_{2b}^3$   
given by \eqref{type2bGH}, \eqref{type2bGRa}, \eqref{type2bGRb}
and \eqref{type2bGRc}, respectively.\\ 
\textit{(i.)} Suppose that the Gauss maps of $X_{2b}$ and $R_{2b}^{1}$ are same. 
From the equations \eqref{type2bGH} and \eqref{type2bGRa}, we get $z(u)=0$ or $z'(u)=0$ which implies $\nu_{R_{2b}^1}=0$.
That is a contradiction. Hence, their Gauss maps are definitely different.
\\
\textit{(ii.)} Suppose that the surfaces $X_{2b}$ and $R_{2b}^{2}$ have the same Gauss maps. 
Comparing \eqref{type2bGH} and \eqref{type2bGRb}, we find $z(u)=0$ or $z'(u)=0$ which can't be possible. 
If $z(u)\neq0$ or $z'(u)\neq0$, then we have the following system of equations:
\begin{align}
\label{type2beq12}
    \lambda y'&=0,\\
\label{type2beq13}
    x'z\sinh{v}-\lambda z'\cosh{v}&=-azz'\cosh\left(v+\int\frac{\lambda x'}{\lambda^2-z^2}du\right), \\
\label{type2beq14}
    x'z\cosh{v}-\lambda z'\sinh{v}&=-azz'\sinh\left(v+\int\frac{\lambda x'}{\lambda^2-z^2}du\right),\\
\label{type2beq15}
    y'z\sinh{v}&=-bzz'\cosh\left(v+\int\frac{\lambda x'}{\lambda^2-z^2}du\right),\\
\label{type2beq16}
     y'z\cosh{v}&=-bzz'\sinh\left(v+\int\frac{\lambda x'}{\lambda^2-z^2}du\right).
\end{align}
Due to $\lambda \neq 0 $,  the equation \eqref{type2beq12} gives $y'(u)=0$.
Then, from the equations \eqref{type2beq15} and \eqref{type2beq16}
imply $b(u)=0$.
Therefore, it can be easily seen that 
the surfaces $X_{2b}$ and $R_{2b}^{2}$ are hyperplanar,
that is, they are lying in $\mathbb{E}^3_1$. 
Moreover, the equations \eqref{type2beq5} and \eqref{type2bizm1} imply that $H_{1}^{X_{2b}}=H_{1}^{R_{2b}^2}=0$. Also, from the equations \eqref{type2beq5} and \eqref{type2bizm1}, we have
\begin{align}
\label{type2beq18}
\begin{split}
H_{2}^{X_{2b}}&=\frac{z^2(z'(x'z'+zx'')+x'(x'^2-zz''))-
\lambda^2(z'(2x'z'+zx'')-x'zz'')}{2{(z^2(x'^2+z'^2)-\lambda^2z'^2)}^{3/2}},\\
H_{2}^{R_{2b}^2}&=
\frac{a'(z^2-\lambda^2)-azz'(a^2-1)}
{2zz'\sqrt{(\lambda^2-z^2)(1-a^2)^3}}.
\end{split}
\end{align}
Using $b(u)=0$, from the equation \eqref{type2eqR2a} we have
\begin{equation}
 \label{type2beq18a}   
 a^2(u)=\frac{\lambda^2 z'^2(u)z-z^2(u)x'^2(u)}{z^2(u)z'^2(u)}.
\end{equation}
Using the equation \eqref{type2beq18a} in \eqref{type2beq18}, we get
\begin{equation}
 \label{type2beq18b}
H_{2}^{R_{2b}^{2}}=\frac{x'z'^2(2\lambda^2-z^2)
-z^2x'^3+z(\lambda^2-z^2)(z'x''-x'z''))}{2(z^2(x'^2+z'^2)-\lambda^2z'^2)^{3/2}\sqrt{(\lambda^2z'^2-z^2x'^2)(\lambda^2-z^2)}} 
\end{equation}
which implies $H_2^{R_{2b}^{2}}=-\frac{x'z^2}{\sqrt{(\lambda^2z'^2-z^2x'^2)(\lambda^2-z^2)}}H_2^{X_{2b}}$. 
Moreover, using equations \eqref{type2beq13} and \eqref{type2beq14}, we obtain the following equations
\begin{align}
\label{type2beq21}
    -x'z=azz'\sinh{\left(\int\frac{\lambda x'}{\lambda^2-z^2}du\right)},\\
\label{type2beq22}
    \lambda z'=azz'\cosh{\left(\int\frac{\lambda x'}{\lambda^2-z^2}du\right)}.    
\end{align}
Considering the equations \eqref{type2beq21} and \eqref{type2beq22} together, we have
\begin{equation}
\label{type2beq23}
    -\frac{\lambda z'}
    {x'z}=\coth{\left(\int\frac{\lambda x'}{\lambda^2-z^2}du\right)}.
\end{equation}
If we take the derivative of the equation \eqref{type2beq22} with respect to $u$, the equation \eqref{type2beq22} becomes
\begin{equation}
\label{type2beq1d}
    x'z'^2(2\lambda^2-z^2)
-z^2x'^3+z(\lambda^2-z^2)(z'x''-x'z''))=0
\end{equation}
which implies $H_{2}^{X_{2b}}=H_{2}^{R_{2b}^{2}}=0$ in the equation 
\eqref{type2beq1d}. 
Thus, we get the desired results.
Since $R_{2b}^{2}$ is minimal, from the equation \eqref{type2beq1d} 
we have the following differential equation
\begin{equation}
\label{coreqtype2b}
(z^2-\lambda^2) a'+zz'a=zz'a^3
\end{equation}
which is a Bernoulli equation. 
Then, the general solution of this equation is found as 
\begin{equation}
\label{cortype2beq1}
   a^2=\frac{1}{1+c_3(\lambda^2-z^2)}   
\end{equation}
for an arbitrary positive constant $c_3$. Comparing the equations \eqref{type2beq18a} and \eqref{cortype2beq1}, we get
\begin{equation}
 x(u)
 =\pm\sqrt{1+c_3\lambda^2}
 \int{\frac{z'(u)}{z(u)}
 \sqrt{\frac{\lambda^2-z^2(u)}{1+c_3(\lambda^2-z^2)}}}du.  
\end{equation}
whose solution is given by \eqref{eqcor2c} for $c_3>0$. 
Moreover, using the last component of $R_{2b}^2(u,v)$ in \eqref{type2bizaGR}, we have 
\begin{equation}
\int\frac{z(u)z'(u)}{\sqrt{(\lambda^2-z^2(u))(1+c_3(\lambda^2-z^2(u)))}}du
=\pm\frac{1}{\sqrt{c_3}}\arcsinh{\left(\sqrt{c_3(\lambda^2-z^2(u))}\right)}+c_4
\end{equation}
for any arbitrary constant $c_4$.
\\
\textit{(iii.)}
Suppose that the surfaces $X_{2b}$ and $R_{2b}^{3}$ have the same Gauss maps. 
From \eqref{type2bGH} and \eqref{type2bGRc}, we get the following system of equations:
\begin{align}
\label{type2ceq12}
    \lambda y'&=0,\\
\label{type2ceq13}
    x'z\sinh{v}-\lambda z'\cosh{v}&=azz'\sinh\left(v-\int\frac{\lambda x'}{z^2-\lambda^2}du\right), \\
\label{type2ceq14}
    x'z\cosh{v}-\lambda z'\sinh{v}&=azz'\cosh\left(v-\int\frac{\lambda x'}{z^2-\lambda^2}du\right),\\
\label{type2ceq15}
    y'z\sinh{v}&=bzz'\sinh\left(v-\int\frac{\lambda x'}{z^2-\lambda^2}du\right),\\
\label{type2ceq16}
     y'z\cosh{v}&=bzz'\cosh\left(v-\int\frac{\lambda x'}{z^2-\lambda^2}du\right).
\end{align}
Due to $\lambda \neq 0 $,  the equation \eqref{type2ceq12} gives $y'(u)=0$.
Then, from the equations \eqref{type2ceq15} and \eqref{type2ceq16}
imply $b(u)=0$.
Therefore, it can be easily seen that 
the surfaces $X_{2b}$ and $R_{2b}^{3}$ are hyperplanar,
that is, they are lying in $\mathbb{E}^3_1$. 
Moreover, the equations \eqref{type2beq5} and \eqref{type2bizm2} imply 
that $H_{1}^{X_{2b}}=H_{1}^{R_{2b}^3}=0$. Also, from the equations \eqref{type2beq5} and \eqref{type2bizm2}, 
we have
\begin{align}
\label{type2ceq18}
\begin{split}
H_{2}^{X_{2b}}&=\frac{z^2(z'(x'z'+zx'')+x'(x'^2-zz''))-
\lambda^2(z'(2x'z'+zx'')-x'zz'')}{2(z^2(x'^2+z'^2)-\lambda^2z'^2)^{3/2}},\\
H_{2}^{R_{2b}^3}&=
\frac{a'(z^2-\lambda^2)+azz'(1+a^2)}
{2zz'\sqrt{(z^2-\lambda^2)(1+a^2)^3}}.
\end{split}
\end{align}
Using $b(u)=0$, from the equation \eqref{type2eqR3a} we have
\begin{equation}
 \label{type2ceq18a}   
 a^2(u)=\frac{z^2(u)x'^2(u)-\lambda^2 z'^2(u)}{z^2(u)z'^2(u)}.
\end{equation}
Using the equation \eqref{type2ceq18a} in \eqref{type2ceq18}, we get
\begin{equation}
 \label{type2ceq18b}
H_{2}^{R_{2b}^{3}}=\frac{z^2x'(z(z^2-\lambda^2)(x'z''-z'x'')
+x'z'^2(2\lambda^2-z^2)-z^2x'^3)}{2(z^2(x'^2+z'^2)-\lambda^2z'^2)^{3/2}\sqrt{(z^2-\lambda^2)(z^2x'^2-\lambda^2z'^2)}}    
\end{equation}
which implies $H_2^{R_{2b}^{3}}=-\frac{z^2x'}{\sqrt{(z^2x'^2-\lambda^2z'^2)(z^2-\lambda^2)}}H_2^{X_{2b}}$. 
Moreover, using equations \eqref{type2ceq13} and \eqref{type2ceq14}, we obtain the following equations
\textbf{\begin{align}
\label{type2ceq21}
    x'z=azz'\cosh{\left(\int\frac{\lambda x'}{z^2-\lambda^2}du\right)},\\
\label{type2ceq22}
    \lambda z'=azz'\sinh{\left(\int\frac{\lambda x'}{z^2-\lambda^2}du\right)}.    
\end{align}}
Considering the equations \eqref{type2ceq21} and \eqref{type2ceq22} together, we have
\begin{equation}
\label{type2ceq23}
    \frac{x'z}
    {\lambda z'}=\coth{\left(\int\frac{\lambda x'}{z^2-\lambda^2}du\right)}.
\end{equation}
If we take the derivative of the equation \eqref{type2ceq23} with respect to $u$, \eqref{type2ceq23} becomes
\begin{equation}
\label{type2ceq24}
    x'z'^2(z^2-2\lambda^2)
+z^2x'^3+z(\lambda^2-z^2)(x'z''-x''z')=0
\end{equation}
which implies $H_{2}^{X_{2b}}=H_{2}^{R_{2b}^{3}}=0$ in the equation 
\eqref{type2ceq18}. 
Thus, we get the desired results.
Since $R_{2b}^{3}$ is minimal, from the equation \eqref{type2ceq18} 
we have the following differential equation
\begin{equation}
\label{coreqtype2c}
(z^2-\lambda^2)a'+zz'a=-zz'a^3
\end{equation}
which is a Bernoulli equation. 
Then, the general solution of this equation is found as 
\begin{equation}
\label{cortype2ceq1}
   a^2=\frac{1}{c_3(z^2-\lambda^2)-1}   
\end{equation}
for an arbitrary positive constant $c_3$. 
Comparing the equations \eqref{type2ceq18a} and \eqref{cortype2ceq1}, we get
\begin{equation}
 x(u)
 =\pm\sqrt{1+c_3\lambda^2}
 \int{\frac{z'(u)}{z(u)}
 \sqrt{\frac{z^2(u)-\lambda^2}{c_3(z^2-\lambda^2)-1}}}du.  
\end{equation}
whose solution is given by \eqref{eqcor2c} for $c_3>0$. 
Moreover, using the last component of $R_{2b}^3(u,v)$ in \eqref{type2bizaGR}, we have 
\begin{equation}
\int\frac{z(u)z'(u)}{\sqrt{(z^2(u)-\lambda^2)(c_3(z^2(u)-\lambda^2)-1)}}du
=\pm\frac{1}{\sqrt{c_3}}\arccosh{\left(\sqrt{c_3(z^2(u)-\lambda^2)}\right)}+c_4
\end{equation}
for any arbitrary constant $c_4$.
\end{proof}

\begin{remark}
If $x'(u)=0$ for $u\in I\subset\mathbb{R}$, 
then the timelike helicoidal surface given by \eqref{eq12} reduces to the timelike right helicoidal surface in $\mathbb{E}^4_1$. 
On the other hand, $W=(\lambda^2-z^2(u))(y'^2(u)+z'^2(u))<0$ for $\lambda^2-z^2(u)<0$.
Thus, from Theorem \ref{type2bBour}, 
we get the timelike rotational surface $R_{2b}^{3}(u,v)$
which are isometric to the timelike right helicoidal surface in $\mathbb{E}^4_1$. 
Also, Theorem \ref{type2bthrm6} implies that if the timelike right helicoidal and $R_{2b}^{3}$
have the same Gauss map, then we get $a^2(u)=-\frac{\lambda^2}{z^2(u)}$ which gives a contradiction.
Thus, they have the different Gauss maps. 
\end{remark}

Now, we give an example by using Theorem \ref{type2bthrm6}. 

\begin{example}
If we choose $z(u)=u$, $c_3=\frac{1}{2}$, $\lambda=1$ and $c_4=0$, 
then isometric surfaces in \eqref{type2bizaGH} and \eqref{type2bizaGR} 
are given as follows
\begin{equation*}
\label{type2bizaGHn}
X_{2b}(u,v) =
\left(\frac{\sqrt{3}}{2}\arcsinh{\sqrt{\frac{u^2-3}{2}}}
-\frac{1}{2}\arctanh{\sqrt{\frac{u^2-3}{3u^2-3}}}+ v,u\cosh v,u\sinh v\right)
\end{equation*}
and 
\begin{equation*}
\label{type2bizaGRn}
  R_{2b}^{3}(u,v)=
\left(
\begin{array}{c}
\sqrt{2}\arccosh{\sqrt{\frac{u^2-1}{2}}}\\
\sqrt{u^2-1}\cosh{\left(v-\arctanh{\sqrt{\frac{u^2-3}{3u^2-3}}}\right)}\\
\sqrt{u^2-1}\sinh{\left(v-\arctanh{\sqrt{\frac{u^2-3}{3u^2-3}}}\right)}
\end{array}
\right). 
\end{equation*}
For $2\leq u\leq 8$ and $-1\leq v\leq1$,  
the graphs of timelike helicoidal surface $X_{2b}$ and timelike rotational surface $R_{2b}^3$ 
in $\mathbb{E}^3_1$ can be plotted by using Mathematica 10.4 as follows:
\begin{figure}[h]
\centering
\begin{subfigure}{.35\textwidth}
  \centering
  \includegraphics[width=6cm, height=7cm]{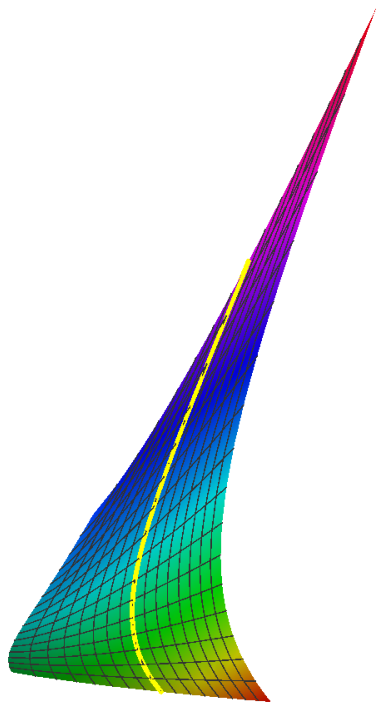}
  \caption{}
  \label{fig:sub1}
\end{subfigure}%
\begin{subfigure}{.35\textwidth}
  \centering
  \includegraphics[width=6cm, height=7cm]{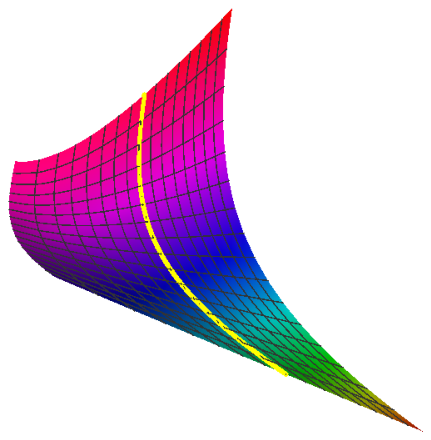}
  \caption{}
  \label{fig:sub2}
\end{subfigure}
\caption{(A) Timelike helicoidal surface of type IIb; timelike helix and (B)  Timelike rotational surface; timelike hyperbola.}
\end{figure}
\end{example}
\section{Helicoidal Surface of Type III}
Let $\{\eta_{1},\eta_{2},\mathbf{\xi}_{3},\mathbf{\xi}_{4}\}$ be the pseudo--orthonormal basis
of $\mathbb{E}^4_{1}$ such that $\mathbf{\xi}_{3}=\frac{1}{\sqrt{2}}(\eta_{4}-\eta_{3})$ 
and $\mathbf{\xi}_{4}=\frac{1}{\sqrt{2}}(\eta_{3}+\eta_{4})$. 
We choose as a lightlike $2-$plane $P_{3}=span\{\eta_{1},\mathbf{\xi}_{3}\}$, 
a hyperplane $\Pi_{3}=span\{\eta_{1},\mathbf{\xi}_{3},\mathbf{\xi}_{4}\}$ 
and a line $l_{3}=span\{\mathbf{\xi}_{3}\}$. 
Then, the orthogonal transformation $T_{3}$ of $\mathbb{E}^4_{1}$ 
which leaves the lightlike plane $P_{3}$ invariant is given by 
$T_{3}(\eta_{1})=\eta_{1},\ T_{3}(\eta_{2})=\eta_{2}+\sqrt{2}v\mathbf{\xi}_{3},\ T_{3}(\mathbf{\xi}_{3})=\mathbf{\xi}_{3}$ and $T_{3}(\mathbf{\xi}_{4})=\sqrt{2}v\eta_{2}+v^{2}\mathbf{\xi}_{3}+\mathbf{\xi}_{4}$. 
We suppose that $\beta_{3}(u)=x(u)\eta_{1}+z(u)\mathbf{\xi}_{3}+w(u)\mathbf{\xi}_{4}$ is a regular curve, where $w(u)\neq 0$. Thus, the parametrization of $X_{3}$ (called as the helicoidal surface of type III) which is obtained a rotation 
of the curve $\beta_{3}$ which leaves the lightlike plane ${P}_{3}$ 
pointwise fixed followed by the translation along $l_{3}$ as follows:
\begin{equation}
\label{eq12}
X_{3}(u,v)=x(u)\eta_{1}+\sqrt{2}vw(u)\eta_{2}+(z(u)+v^{2}w(u)
+\lambda v)\mathbf{\xi}_{3}+w(u)\mathbf{\xi}_{4},
\end{equation}
where, $v\in\mathbb{R}$ and $\lambda\in \mathbb{R^{+}}$. 
When $w$ is a constant function, $X_{3}$ is called as right helicoidal surface of type III.
For $\lambda=0$, the helicoidal surface which is given by \eqref{eq12} reduces 
to the rotational surface of parabolic type in $\mathbb{E}^4_1$ (see \cite{Dursun} and \cite{Bektas}).

By a direct calculation, we get the induced metric of $X_3$ given as follows.
\begin{equation}
\label{type3eq2}
ds^2_{X_3}=(x'^2(u)-2w'(u)z'(u))du^2 -2\lambda w'(u)dudv + 2w^2(u)dv^2.
\end{equation}
Due to the fact that $X_3$ is a timelike helicoidal surface in $\mathbb{E}^4_1$, 
we have ${W}=2w^2(u)(x'^2(u)-2w'(u)z'(u))-\lambda^2w'^2(u)<0$
for all $u\in I\subset\mathbb{R}$. 
Then, we choose an orthonormal frame field 
$\{e_{1},e_{2}, N_{1}, N_{2}\}$ on $X_3$ in $\mathbb{E}^4_1$ 
such that $e_1, e_2$ are tangent to $X_3$ and $N_1, N_2$ are normal to $X_3$ as follows.
\begin{align}
\label{type3eq3}
\begin{split}
e_{1}&=\frac{1}{\sqrt{\epsilon g_{11}}}{{X_3}_u},\;\;\;\;
e_{2}=
\frac{1}{\sqrt{-\epsilon Wg_{11}}}(g_{11}{X_3}_v-{g_{12}{X_3}_u}), \\
N_{1}&=\eta_{1}+\frac{x'}{w'}\xi_{3},\\
N_{2}=&\frac{1}{w'\sqrt{-W}}\bigg(\sqrt{2}x'ww'\eta_{1}+(\lambda w'^2+2vww'^2)\eta_{2}
+\sqrt{2}(\lambda vw'^2+v^2ww'^2+wx'^2\\
&-ww'z')\xi_{3}+\sqrt{2}ww'^2\xi_{4}\bigg),
\end{split}
\end{align}
where $\langle e_{1},e_{1}\rangle=-\langle e_{2},e_{2}\rangle=\epsilon$ and $\langle N_{1},N_{1}\rangle=\langle N_{2},N_{2}\rangle=1$.
It can be easily seen that $X_3$ has a spacelike meridian curve for $\epsilon=1$. 
Otherwise, it has a timelike meridian curve. 
By direct computations, we get the coefficients of the second fundamental form given as follows.
\begin{align}
\label{type3eq4}
    \begin{split}
     &b_{11}^1=\frac{x'' w'- x' w''}{w'},\;\;\; 
     b_{12}^1=b_{21}^1=b_{22}^1=0,\\
     &b_{11}^2
     =\frac{\sqrt{2}w(x'x''w'-x'^2w''+w'(z'w''-w'z''))}{w'\sqrt{-W}},\;\; 
     b_{12}^2=b_{21}^2=\frac{\sqrt{2}\lambda w'^2}{\sqrt{-W}},\;\;
     b_{22}^2=-\frac{2\sqrt{2}w^2w'}{\sqrt{-W}}.
    \end{split}
\end{align}
Thus, we find the components of mean curvature vector $H$ of $X_3$ in $\mathbb{E}^4_1$ as 
\begin{align}
\label{type3eq5} 
\begin{split}
&H_{1}^{X_{3}}=\frac{w^2(x''w'-x'w'')}
{w'W},\\
&H_{2}^{X_{3}}=\frac{-\sqrt{2}(\lambda^2w'^4
+2w^2w'^3z'-w^3x'^2w''+w^3w'(z'w''+x'x'')
-w^2w'^2(x'^2+wz''))}
{w'(-W)^{3/2}}.
\end{split}
\end{align}
We note that if $w'(u)=0$ for $u\in I$, then $W>0$. Thus, $w'(u)$ must be different than zero for $u\in I$.

\subsection{Bour's Theorem and the Gauss map for helicoidal surface of type III}
In this section, we study on Bour's theorem for timelike helicoidal surface of type III in $\mathbb{E}^4_1$ 
and we analyse the Gauss maps of isometric pair of surfaces. 

\begin{theorem}
\label{type3Bour}
A timelike helicoidal surface of type III in $\mathbb{E}^4_{1}$ given by \eqref{eq12} is isometric to one of 
the following timelike rotational surfaces in $\mathbb{E}^4_1$:
\begin{equation}
\label{type3izR}
\begin{split}
R_{3}(u,v)&=
\int{a(u)w'(u)du}\eta_{1}+
\sqrt{2}w(u)\left(v+\frac{\lambda}{2w(u)}\right)\eta_{2}\\
&+\left(\int{b(u)w'(u)du}+w(u)\left(v+\frac{\lambda}{2w(u)}\right)^2\right)\xi_{3}+w(u)\xi_{4}    
\end{split}
\end{equation}
so that spacelike helices on the timelike helicoidal surface of type III correspond to parallel spacelike parabolas 
on the timelike rotational surfaces, 
where $a(u)$ and $b(u)$ are differentiable functions satisfying the following equation:
\begin{equation}
 \label{type3eq7}
    a^2(u)-2b(u)=\frac{x'^2(u)-2w'(u)z'(u)}
    {w'^2(u)}-\frac{\lambda^2}{2w^2(u)}
 \end{equation}
with $w'(u)\not=0$ for all $u\in I\subset\mathbb{R}$.
\end{theorem}

\begin{proof}
Assume that $X_3$ is a timelike helicoidal surface of type III in $\mathbb{E}^4_1$ defined by \eqref{eq12}. 
Then, we have the induced metric of $X_{3}$ given by \eqref{type3eq2}. 
Now, we will find new coordinates $\bar{u},\bar{v}$ such that the metric becomes 
\begin{equation}
    ds^2_{X_3}=F(\bar{u})du^2+G(\bar{u})d\bar{v}^2,
\end{equation}
where $F({\bar{u}})$ and $G(\bar{u})$ are smooth functions. 
Set $\bar{u}=u$ and $\overline{v}=v+\frac{\lambda}{2w(u)}$.
Since Jacobian $\displaystyle{\frac{\partial(\bar{u},\bar{v})}{\partial(u,v)}}$ is nonzero,
it follows that $\{\bar{u}, \bar{v}\}$ are new parameters of $X_3$. 
According to the new parameters, the equation \eqref{type3eq2} becomes
\begin{equation}
\label{type3eq10}
 ds^2 _{X_3}=\left(x'^2(u)-2w'(u)z'(u)-\frac{\lambda^2 w'^2(u)}{2w^2(u)}\right)du^2 + 2w^2(u)d\overline{v}^2.
\end{equation}
On the other hand, 
the timelike rotational surface $R_3$ in $\mathbb{E}^4_{1}$
related to $X_3$ is given by 
\begin{equation}
\label{type3rtk}
R_{3}(k,t)=n(k)\eta_{1}+\sqrt{2}tr(k)\eta_{2}+(s(k)+t^2r(k))\xi_{3}+r(k)\xi_{4}.
\end{equation}
We know that the induced metric of $R_3$ is given by
\begin{equation}
\label{type3eq11}
 ds^2 _{R_3}=(\dot{n}^2(k)-2\dot{r}(k)\dot{s}(k))dk^2 + 2r^2(k)dt^2
\end{equation}
with $\dot{n}^2(k)-2\dot{r}(k)\dot{s}(k)<0$. 
From the equations \eqref{type3eq10} and \eqref{type3eq11}, 
we get an isometry between $X_3$ and $R_{3}$ by taking
$\bar{v}=t$, $r(k)=w(u)$ 
and
\begin{equation}
\label{type3eq12}
 \left(x'^2(u)-2w'(u)z'(u)-\frac{\lambda^2 w'^2(u)}{2w^2(u)}\right)du^2
=(\dot{n}^2(k)-2\dot{r}(k)\dot{s}(k)) dk^2.
\end{equation}
Let define $a(u)=\frac{\dot{n}(k)}{\dot{r}(k)}$ and 
$b(u)=\frac{\dot{s}(k)}{\dot{r}(k)}$.
Using these in the equation \eqref{type3eq12}, we obtain the equation \eqref{type3eq7}.
\begin{equation}
  n=\int{a(u)w'(u)du}\;\;\mbox{and}\;\;s=\int{b(u)w'(u)du}.
\end{equation}
Thus, we get an isometric timelike rotational surface $R_3$ given by \eqref{type3izR}. 
Moreover, 
if we choose a spacelike helix $X_3(u_0,v)$ on $X_3$
for an arbitrary constant $u_{0}$, then 
it corresponds to 
$R_{3}(u_0,v)=
\sqrt{2}w_0\left(v+\frac{\lambda}{2w_0}\right)\eta_{2}
+w_0\left(v+\frac{\lambda}{2w_0}\right)^2\xi_{3}+w_0\xi_{4}$.    
If we take $t=v+\frac{\lambda}{2w_0}$, then it can be rewritten 
$\alpha(t)=\sqrt{2}w_0
\left(0,t,-\frac{t^2}{2},\frac{t^2}{2}\right)
+\frac{1}{\sqrt{2}}\left(0,0,w_0,w_0\right).
$
From Definition \ref{circledef}, it can be seen that $\alpha(t)$ 
is a spacelike parabola lying on the $x_3x_4$--plane.
\end{proof}

\begin{lemma}
Let $X_3$ and $R_3$ be timelike surfaces in $\mathbb{E}^4_1$ given by \eqref{eq12} and \eqref{type3izR}, respectively. 
Then, the Gauss maps of them 
\begin{align}
\label{type3GH}
\nu_{X_{3}}=&\frac{\epsilon}{\sqrt{-W}}\Bigg(\sqrt{2}x'w\eta_1\wedge\eta_2
+x'(\lambda+2vw)\eta_1\wedge\xi_3+ \sqrt{2}(v^2ww'-wz'+\lambda vw')\eta_2\wedge\xi_3\notag\\
&-\sqrt{2}ww'\eta_2\wedge\eta_4-w'(\lambda+2vw)\xi_3\wedge\xi_4\Bigg),\\
\label{type3GR}
\nu_{R_{3}}=&\frac{\epsilon ww'}{\sqrt{-W}}
\Bigg(\sqrt{2}a\eta_1\wedge\eta_2
+2a\left(v+\frac{\lambda}{2w}\right)\eta_1\wedge\xi_3+\sqrt{2}\left(\left(v+\frac{\lambda}{2w}\right)^2-b\right)\eta_2\wedge\xi_3\notag\\
&-\sqrt{2}\eta_2\wedge\xi_4
-2\left(v+\frac{\lambda}{2w}\right)\xi_3\wedge\xi_4\Bigg).
\end{align}
\end{lemma}

\begin{proof}
It follows from a direct calculation. 
\end{proof}

\begin{theorem}
A timelike helicoidal surface of type III and a timelike rotational surface in $\mathbb{E}^4_1$ 
given by \eqref{eq12} and \eqref{type3izR}, respectively have the same Gauss map.
\end{theorem}

\begin{proof}
Assume that the surfaces $X_3$ and $R_3$ have the same Gauss map.
Comparing \eqref{type3GH} and \eqref{type3GR}, we get the following system of equations
\begin{align}
\label{type3eq12a}
x'&=aw',\\
\label{type3eq12b}
x'\left(\lambda+2vw\right)&=2aww'\left(v+\frac{\lambda}{2w}\right),\\
\label{type3eq12c}
wz'&=bww'-\frac{\lambda^2w'}{4w}.
\end{align}
From the equations \eqref{type3eq12a} and \eqref{type3eq12c}, we find $a(u)$ and $b(u)$. 
Using these in \eqref{type3eq7}, we can see that they have the same Gauss map.
\end{proof}

We note that T. Ikawa \cite{Ikawa2} studied Bour's theorem for helicoidal surfaces in $\mathbb{E}^3_1$
with lightlike axis and he showed that they have the same Gauss map.

Now, we give an example by using Theorem \ref{type3Bour}. 

\begin{example}
If we choose $x(u)=a(u)=0$, $w(u)=z(u)=u$ and $\lambda=5$, 
then isometric surfaces in \eqref{eq12} and \eqref{type3izR} are given as follows
\begin{equation*}
\label{type3aGH}
X_{3}(u,v)=\sqrt{2}uv\eta_{2} +\left(u+uv^2+5v\right)\xi_{3}+u\xi_{4}
\end{equation*}
and
\begin{equation*}
    R_{3}(u,v)=\left(\sqrt{2}uv+\frac{5}{\sqrt{2}}\right)\eta_{2}
    +\left(u-\frac{25}{4u}+u\left(v+\frac{5}{2u}\right)^2 \right)\xi_{3}+u\xi_{4}.
\end{equation*}
For $-4\leq u\leq 4$ and $-3\leq v\leq3$,  
the graphs of timelike helicoidal surface $X_{3}$ and 
timelike rotational surface $R_{3}$ in $\mathbb{E}^3_1$ can be plotted by using Mathematica 10.4 as follows:
\begin{figure}[h]
\centering
\begin{subfigure}{.3\textwidth}
  \centering
  \includegraphics[width=4cm, height=6cm]{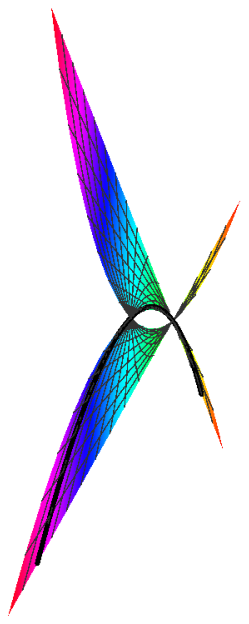}
  \caption{}
  \label{fig:sub1}
\end{subfigure}%
\begin{subfigure}{.3\textwidth}
  \centering
  \includegraphics[width=4cm, height=6cm]{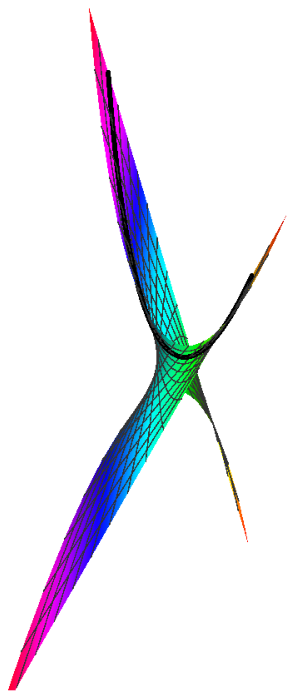}
  \caption{}
  \label{fig:sub2}
\end{subfigure}
\caption{(A) Timelike helicoidal surface of type III; spacelike helix and (B)  Timelike rotational surface; spacelike parabola.}
\end{figure}
\end{example}

\section{Conclusion}

In this paper, we study on Bour's theorem for four kinds of timelike helicoidal surfaces in 4-dimensional Minkowski space. Moreover, we analyse the geometric properties of these isometric surfaces having same Gauss map. Also, we determine the parametrizations of such isometric pair of surfaces. Finally, we give some examples by using Wolfram Mathematica 10.4.

In the future, we will try to determine the helicoidal and rotational surfaces which are isometric according to Bour's theorem whose the mean curvature vectors or their lengths are zero and the Gaussian curvatures are zero, respectively.

\section*{Acknowledgment}
This work is a part of the master thesis of the third author and it is supported by The Scientific and Technological Research Council of Turkey (TUBITAK) under Project 121F211.

\end{document}